\documentclass[11pt,a4paper]{article}
\date{} 
\usepackage{amsfonts,amsmath,amsthm}
\usepackage{amssymb,latexsym}
\usepackage[dvips]{epsfig}
\usepackage{amscd}  
\usepackage[all]{xy}
\usepackage{enumerate}

\newcommand{\ov}{\overline}

\newtheorem{theorem}{Theorem}[section]
\newtheorem{lemma}[theorem]{Lemma}
\newtheorem{cor}[theorem]{Corollary}
\newtheorem{defn}[theorem]{Definition}
\newtheorem{prop}[theorem]{Proposition}
\newtheorem{proposition}[theorem]{Proposition}

\newtheorem{exs}[theorem]{Examples}

\newcommand{\cA}{{\mathcal A}}

\newcommand{\cC}{{\mathcal C}}
\newcommand{\cD}{{\mathcal D}}

\newcommand{\cF}{{\mathcal F}}

\newcommand{\cH}{{\mathcal H}}

\newcommand{\cL}{{\mathcal L}}

\newcommand{\cO}{{\mathcal O}}
\newcommand{\cP}{{\mathcal P}}
\newcommand{\cR}{{\mathcal R}}

\newtheorem{remark}[theorem]{Remark}

\newtheorem{conj}[theorem]{Conjecture}

\newcommand{\mg}{\mathfrak{g}}
\newcommand{\mh}{\mathfrak{h}}
\newcommand{\mb}{\mathfrak{b}}

\newcommand{\mS}{\mathbb{S}}
\newcommand{\mV}{\mathbb{V}}

\newcommand{\mC}{\mathbb{C}}

\newcommand{\mZ}{\mathbb{Z}}

\newcommand{\ie}{{\it i.e. }}
\newcommand{\eg}{{\it e.g. }}

\newcommand{\la}{\lambda}

\newcommand{\HOM}{\operatorname{Hom}}

\newcommand{\EXT}{\operatorname{Ext}}

\newcommand{\END}{\operatorname{End}}

\newcommand{\MOD}{\operatorname{-mod}}
\newcommand{\MOF}{\operatorname{-mod}}

\newcommand{\op}{\operatorname}

\newcommand{\KER}{\operatorname{ker}}

\newcommand{\COKER}{\operatorname{coker}}

\newcommand{\DIM}{\operatorname{dim}}

\newcommand{\opp}{^{\mbox\sl opp}}

\newcommand{\surj}{\mbox{$\rightarrow\!\!\!\!\!\rightarrow$}}

\newcommand{\ID}{\operatorname{ID}}

\newcommand{\p}{\mathfrak{p}}

\numberwithin{equation}{section}

\title{Projective-injective modules, Serre functors and symmetric algebras} 
\author{Volodymyr Mazorchuk and Catharina Stroppel} 

\begin{document}

\maketitle

\baselineskip 14pt

\def\C{\mathbb C}
\def\R{\mathbb R}
\def\N{\mathbb N}
\def\Z{\mathbb Z}
\def\Q{\mathbb Q}
\def\g{\mathfrak g}
\def\p{\mathfrak p}
\def\h{\mathfrak h}
\def\n{\mathfrak n}
\newcommand{\CA}{\operatorname{Coapp}}
\newcommand{\Ext}{\operatorname{Ext}}
\newcommand{\End}{\operatorname{End}}
\newcommand{\add}{\operatorname{add}}
\newcommand{\KZ}{Knizhnik-Zamolodchikov}

\def\l{\lbrace}
\def\r{\rbrace}
\def\lra{\longrightarrow}
\def\Hom{\mathrm{Hom}}

\def\mc{\mathcal} 
\def\mf{\mathfrak}
\def\mb{\mathbb}
\def\cO{\mc{O}}   
\def\cC{\mc{C}}
\def\sln{\mathfrak{sl}(n)} 

\begin{abstract} 
We describe Serre functors for (generalisations of) the category $\cO$ associated with
a semi-simple complex Lie algebra. In our approach, projective-injective
modules, that is modules which are both, projective and injective, play an important role. They control the Serre functor in the case of a 
quasi-hereditary algebra having a double centraliser given by a
projective-injective module whose endomorphism ring is a symmetric algebra. As an application of the double centraliser property
together with our description of Serre functors, we prove three conjectures of Khovanov about the 
projective-injective modules in the parabolic category $\mc{O}_0^{\mu}(\mathfrak{sl}_n)$. 
\end{abstract} 

\tableofcontents
\section{Introduction}

Symmetric algebras are particularly well-behaved algebras with several special
properties. We first want to recall two of these properties and then discuss
to which extend they can be generalised and transferred to other finite
dimensional algebras. 
A finite dimensional algebra $A$ is called {\it symmetric} if there is an
$A$-bimodule isomorphism $A\cong A^*$. In particular, projective $A$-modules
are always injective. On the other hand, since $A\cong A^*$ as $A$-bimodules, 
the isomorphism of vector spaces 
\begin{equation*}
\HOM_A(A,A)\cong\HOM_A(A, A^*\otimes_A A)\cong\HOM_A(A,A)^*
\end{equation*}
gives rise to an isomorphism 
\begin{equation*}
\HOM_A(M,N)\cong\HOM_A(M, A^*\otimes_A N)\cong\HOM_A(N,M)^*
\end{equation*} 
for any free $A$-modules $M$ and $N$. Moreover, this isomorphism is functorial
in $M$ and $N$. \\
The question we want to ask now is whether these two properties can be
transferred somehow to a more general class of algebras. We will consider finite dimensional algebras which have a symmetric subalgebra which, in some sense, can control the representation
theory of the original algebra. If we are lucky, $A$ has ``enough''
projective modules which are also injective. We fix a system of
representatives for the isomorphism classes of indecomposable
projective-injective (\ie at the same time projective {\it and} injective)
$A$-modules. Then, instead of considering $A$
itself we propose to consider the direct sum, say $Q$, of all modules from
this fixed system. The endomorphism ring $\END_A(Q)$ is at least a Frobenius algebra. The following questions naturally
arise
\begin{itemize}
\item Is $\END_A(Q)$ a symmetric algebra?
\item How much information about the category of finitely generated
  $A$-modules is already encoded in $Q$ and $\END_A(Q)$?
\item Is there an isomorphism, functorial in both entries,  
  \begin{equation*}
    \HOM_A(P_1, P_2)\cong\HOM_A(P_1,A^*\otimes_A P_2)\cong \HOM_A(P_2,P_1)^*
  \end{equation*}
for any projective-injective $A$-modules $P_1$ and $P_2$?
\end{itemize}

In general, the first question is very difficult to answer. Concerning the
second question, we will describe the situation where all the information
about $A\MOD$ is already contained in $Q$ and $\END_A(Q)$. This is given by
the so-called double centraliser property which we will explain more precisely
shortly. The last question relates to the existence and description of a so-called Serre
functor, motivated by Serre's duality for sheaves on projective
varieties. The purpose of this paper is to answer the above three questions for
certain algebras appearing in representation theory.

To substantiate and to specify our approach we would like to recall the important role
projective-injective modules play in representation theory, in particular in different versions of
the category $\cO$. The striking example is the case of an integral block of
the category $\cO$ for a semisimple complex Lie algebra. By Soergel's result
(\cite[Endomorphismensatz and Struktursatz]{Sperv}) such a block can be
completely described by (a certain subcategory of) the category of finitely generated modules over the
endomorphism  ring of the (unique up to isomorphism) indecomposable
projective-injective module in this block. This idea was generalised and
formalised in \cite{KSX} as the so-called {\em double centraliser property}. In this language, Soergel's result 
could be stated as follows: Let $A$ be a finite dimensional algebra such that 
$A\MOD$ is equivalent to the block of the category $\cO$ in question. Then for 
the indecomposable projective-injective $A$-module $Q$ we have
$A=\END_{\END_A(Q)}(Q)$. Note that \cite[Endomorphismensatz]{Sperv} implies
that $\END_A(Q)$ is a symmetric algebra. 

Using the Ringel duality functor it is easy to see (Corollary~\ref{doublecentralizer}) 
that there is {\it always} a tilting module $T$ having the double centraliser
property above. This module $T$ need not to be projective-injective in general.
Nevertheless, there are many known examples where we 
have the particularly handy situation of the double centraliser property with respect to a
projective-injective tilting module $T$. We recall such examples in Section~\ref{section1}. 
Since in these cases the category $A\MOD$ is completely determined by $\END_A(T)$ and $T$, 
it follows directly that the centre of $A\MOD$ is isomorphic to the
centre of the endomorphism ring of $T$ (in particular
\cite[Conjecture~4]{Kh1} follows, see Theorem~\ref{t5}). 

Motivated by Serre's duality, there is the notion of a Serre functor for any 
$\Bbbk$-linear category with finite dimensional homomorphism spaces (see 
Section~\ref{section2}). Kapranov and Bondal (\cite{BoKa}) showed that 
the bounded derived category $\cD^b(A)$ for any finite dimensional algebra $A$
of finite global dimension admits a Serre functor. In fact, the existence of a 
Serre functor is equivalent to the finiteness of the  global dimension of $A$ 
and also to the existence of Auslander-Reiten triangles  (\cite{RvB}, \cite{Happeltri}). 
It is well-known that in the latter case the Serre functor is the  
left derived of the Nakayama functor 
(see \eg \cite[page 37]{Happel}), that is of the functor isomorphic to 
$A^\ast\otimes_A\bullet$. However, if the algebra $A$ is not explicitly given,
the Serre functor might be hard to compute.

Nevertheless, in some cases the Serre functor for $\cD^b(A)$ can be explicitly described, using
for instance geometric or functorial methods. For example, in \cite{BoKa} it was 
conjectured that the Serre functor of the bounded derived category of perverse sheaves on flag 
varieties is given by a geometrically defined intertwining functor. 
This was our motivation to study the Serre functor of the bounded derived
category of (integral blocks of) $\cO$, associated 
with the corresponding semi-simple Lie algebra, which is equivalent to the category of 
perverse sheaves in question. 
The original conjecture has recently been proved by Beilinson, Bezrukavnikov and 
Mirkovic in \cite{BBM}. In the present paper we explicitly construct the Serre functor 
for the bounded derived category of any integral block of $\cO$ using the twisting functors, defined in \cite{Arkhipov} and
studied \eg in \cite{AS}. Our approach is purely algebraic and does not require
the explicit knowledge of the associative algebra associated to $\cO$. As a (very unexpected) consequence we obtain an isomorphism between a certain
composition of twisting functors and a certain composition  of Irving's shuffling 
functors (see Corollary~\ref{Arkhipov}).

We further apply this result to construct the Serre functors of the bounded
derived categories of (integral blocks of) of any parabolic
category $\cO^{\mathfrak{p}}$ in the sense of Rocha-Caridi, \cite{RC}.
 Using the explicit description of the Serre functor (in terms of shuffling
 functors) we prove in Theorem~\ref{symmetry} that
the endomorphism algebra of the sum of all indecomposable projective-injective
modules in $\cO^\p$ is symmetric, which was conjectured by Khovanov. 
One of the motivations to consider the category of
projective-injective modules in $\cO^{\p}$ is to find a precise
connection between Khovanov's categorification of the Jones polynomial 
(\cite{Khovanov}) and the categorification of the Jones polynomial via representation 
theory of the Lie algebra $\mathfrak{sl}_n$ (as proposed in \cite{BFK} and
proved in \cite{Sduke}). It might be possible to simplify the approach in
\cite{BFK} and \cite{Sduke} by working with these symmetric endomorphism
algebras. Moreover, from a topological point of view it seems to be much more
natural and plausible to work with symmetric (or at least Frobenius) 
algebras to construct knot invariants instead of the complicated algebras
which describe the integral blocks of the (parabolic) category $\cO$.  

We expect that the Serre
functor for the category $\cO$ for rational Cherednik algebras can be
constructed in a similar way via twisting functors as the Serre functor for the
Bernstein-Gelfand-Gelfand category $\cO$. However, we
are not able to prove this, mainly because of the lack of translation
functors. Nevertheless, we give a description of the Serre functor for the category 
$\cO$ for rational Cherednik algebras via partial coapproximation with respect
to the direct sum of all indecomposable projective-injective modules (see \cite[2.5]{KM}). 
The proof however uses the
fact that Hecke algebras are symmetric and the properties of the
\KZ-functor. Using \cite[Remark 5.10]{GGOR} it might be possible to construct
the Serre functor in a different way, which would imply a conceptual proof of
the fact (\cite{CIK}) that the Hecke algebras occurring here are symmetric (see
Conjecture~\ref{conjCheredsym}). 

In the last section of the paper we consider the special case of a parabolic 
category, $\cO_0^{\p}(\mathfrak{sl}_n)$, for the Lie algebra $\mathfrak{sl}_n$. 
For this category we give an easier proof of the main result 
of \cite{Irself} (Theorem~\ref{t7}). As a consequence we show that there is
a double centraliser property with respect to a basic projective-injective module.
This implies \cite[Conjecture~4]{Kh1}.  The parabolic subalgebra $\p$ of
$\mathfrak{sl}_n$ is determined by some composition of $n$. In \cite[Proposition~4.3]{Irself} 
it is shown that indecomposable projective-injective modules in $\cO_0(\mathfrak{sl}_n)^{\p}$
are indexed by the elements of some left cell. The Kazhdan-Lusztig combinatorics of translation
functors, applied to these projective-injective modules, suggests a connection
with Specht modules for the symmetric group $S_n$. It is well known that the
Specht modules which correspond to different compositions of $n$, but giving 
rise to the same partition of $n$, are isomorphic. This observation might have
led M. Khovanov to the conjecture that the endomorphism 
algebras of the basic projective-injective modules in different $\cO_0(\mathfrak{sl}_n)^{\p}$, 
corresponding to the same partition of $n$, are isomorphic. We finish the
paper by proving this conjecture (Theorem~\ref{tparab}).\\

\noindent{\bf Some guidance for the reader:} Section~\ref{section1} recalls
some facts and results on double centralisers for module categories over
standardly stratified algebras. In principle, the content is not new, the
viewpoint might be slightly more general than usual. We formulated
the setup as generally as possible, since we believe that our approach can be
applied to a much wider class of algebras than the ones actually appearing in
the paper. For the reader, however, it might be more approachable to have first a look at the
Examples~\ref{examples}, skip the details of Section~\ref{section1} and focus
on the main result characterising Serre functors which can be found in
Section~\ref{section2}. Applications and concrete descriptions of Serre
functors are given in Section~\ref{secappl}. The deepest result might be
Theorem~\ref{symmetry} stating that the endomorphism ring of a basic
projective-injective module in the principal block of a (parabolic) category
$\cO$ is not only Frobenius, but symmetric.    

\bigskip
\noindent
{\bf Acknowledgements.}
We thank Mikhail Khovanov for sharing his ideas with us
and for many helpful discussions and remarks. We also would like to thank 
Iain Gordon for many useful discussions and explanations and Ken Brown for
useful remarks. We also thank Joshua
Sussan for comments on a preliminary version of the paper. 
The first author was partially supported by The Royal Swedish 
Academy of Sciences, The Swedish Research Council,
and STINT.  The second author was supported by EPSRC.

\section{Standardly stratified structure and double centralisers}
\label{section1}

In representation theory double centraliser properties play an important role. The aim of this section is to recollect known
results from the literature, to emphasise the universal principle behind it, and
to show the significance of projective-injective modules.

Let $\Bbbk$ be an algebraically closed field. Let $A$ be a unital finite dimensional
associative $\Bbbk$-algebra. We denote by $A\MOF$ ($\operatorname{mod}-A$) 
the category of finite dimensional left (resp. right) $A$-modules. In the following we 
will mainly work with left $A$-modules, hence an ``$A$-module'' is always meant to be a
{\it left} $A$-module. For $M$, $N\in A\MOF$ we denote by $\op{Tr}_MN$
the trace of $M$ in $N$ (which is by definition the submodule of $N$,
generated by the images of all morphisms from $M$ to $N$).  

Let $\{L(\la)\}_{\la\in\Lambda}$ be a complete set of representatives for the 
isomorphism classes of simple $A$-modules. For
a simple $A$-module, $L(\la)$, we 
denote by $P(\la)$ its projective cover, and by $I(\la)$ its injective
hull. We assume that there is a partial pre-order $\preceq$  (\ie
a reflexive and transitive binary relation) on $\Lambda$, which we fix.  
Let $P^{\succ\la}=\bigoplus_{\mu\succ\la} P(\mu)$ and 
$P^{\succeq\la}=\bigoplus_{\mu\succeq\la} P(\mu)$. 

With respect to $\preceq$
we define the so-called {\it standard module} $\Delta(\la)$ to be the 
largest quotient  of $P(\la)$ containing only
composition factors of the form $L(\mu)$, where $\la\not\prec\mu$, \ie\ 
$\Delta(\la)=P(\la)/\op{Tr}_{P^{\succ\la}}P(\la)$. 
We also have a {\it proper standard module} $\overline\Delta(\la)$ which is the
largest quotient of $P(\la)$ such that its radical contains only
composition factors of the form $L(\mu)$, where $\la\not\preceq\mu$,  \ie
$\overline{\Delta}(\la)=P(\la)/\op{Tr}_{P^{\succeq\la}}\op{rad}P(\la)$.
Dually, we have the {\it costandard module} $\nabla(\la)$ and the {\it proper
costandard module} $\overline{\nabla}(\la)$.    

We denote by $\cF(\Delta^A)=\cF(\Delta)$ the full subcategory of $A\MOF$ given by all modules
having a filtration, with all subquotients of this filtration being isomorphic to
$\Delta(\la)$ for various $\la\in\Lambda$. 
If $M\in\cF(\Delta)$ then we say that {\it $M$ has a standard flag}. Similarly, we define
$\cF(\overline{\Delta})$, $\cF(\nabla)$, $\cF(\ov\nabla)$, the categories of
modules having a proper standard, a costandard, and a proper costandard flag
respectively.

Let $A$ be a finite dimensional {\it standardly stratified
algebra} as defined in \cite{CPS2}, that is 
\begin{itemize}
\item the kernel of the canonical surjection $P(\la)\surj\Delta(\la)$
has a standard flag; 
\item the kernel of the canonical surjection $\Delta(\la)\surj L(\la)$ has a
filtration with subquotients $L(\mu)$, where $\mu\preceq\la$. 
\end{itemize}

In particular, if $\preceq$ is a partial order and $\Delta(\la)=\overline\Delta(\la)$
for $\la$, then $A$ is {\it quasi-hereditary} (see \cite{CPS,DR1}). If $\preceq$ is a partial 
order and any $\Delta(\la)$ has a proper standard flag, then $A$ is {\it properly stratified} (see \cite{D1}).

We call a module $M$ {\it basic} with respect to some property $\cP$, if $M$ is the
direct sum of pairwise non-isomorphic indecomposable modules with property
$\cP$ and any indecomposable module having this property is isomorphic to a
summand in $M$. For example, a basic projective module in $A\MOF$ is a minimal
projective generator. If $N=\oplus_{i=1}^k N_i^{m_i} $, where
$m_i\in\{1,2,\dots\}$ for all $i$, with $N_i$ indecomposable and pairwise 
non-isomorphic, we set $N_{basic}=\oplus_{i=1}^k N_i$.

\subsection{Tilting modules and Ringel duality}
A {\it tilting module} is an object in $\cF(\Delta)\cap\cF(\ov\nabla)$, and
a {\it cotilting module} is an object in $\cF(\ov\Delta)\cap\cF(\nabla)$.
In \cite{Frisktilting} it is shown that for a standardly stratified
algebra the category $\cF(\Delta)\cap\cF(\ov\nabla)$ is closed under taking
direct summands and that the indecomposable modules in this category are in 
natural bijection with standard modules. Let $T(\la)$ denote the unique 
indecomposable tilting module having a standard flag, where $\Delta(\la)$ occurs 
as a submodule. Let $T=\oplus _{\la\in\Lambda} T(\la)$ be the {\em characteristic tilting
module}. There is the dual notion of {\it cotilting modules}. In general, cotilting modules cannot be classified in the same way
as tilting modules. However, this can be done in the case when the opposite
algebra $A^{opp}$ is also standardly stratified (with respect to the same
partial pre-order), see \cite[4.2]{Frisktilting}. 
For quasi-hereditary algebras cotilting and tilting modules
obviously coincide, but in general they do not have to. The Ringel 
duality functor (as introduced in \cite{Ringel}) was studied in the more
general setup of various stratified algebras for example
in \cite{AHLU} and \cite{Frisktilting}. We will need the following slight
variation of these results: 

\begin{proposition}
\label{Ringeldual}
Let $A$ be a standardly stratified algebra. Then the Ringel dual 
$R(A)=\END_A(T)$ is standardly stratified and the contravariant 
functor $\mathtt{R}=\HOM_A(\bullet,T):A\MOF\rightarrow \op{mod-}R(A)$
satisfies the following properties:
\begin{enumerate}[(1)]
\item $\mathtt{R}$ maps tilting modules to projective modules.
\item $\mathtt{R}$ maps projective modules to tilting modules.
\item $\mathtt{R}$ defines an equivalence of categories
$\cF(\Delta^A)\cong\cF(\Delta^{R(A)})$.  
\end{enumerate}
\end{proposition}

\begin{proof}
That the algebra $R(A)$ is standardly stratified follows for example
from \cite[Theorem 5 (iii)]{Frisktilting}. Obviously, $T$ is mapped to 
$\End_A(T)$, hence it is projective. Taking direct summands implies 
the first statement. The last statement is proved analogously to  
\cite[Theorem 2.6 (iv)]{AHLU} (note that the duality $\op{D}$ used there
swaps standard and costandard modules). To prove the second statement let now 
$Q$ be projective, then $\mathtt{R}Q$ has a standard flag. Of course, 
$\EXT^1_A(Q,\Delta(\la))=0$ for any $\la$. Using the last part of the 
proposition we get $\EXT^1_{R(A)}(\mathtt{R}\Delta(\la),\mathtt{R} Q)=0$, even
$\EXT^1_{R(A)}(\Delta^{R(A)}(\la),\mathtt{R} Q)=0$ for any standard module
$\Delta^{R(A)}(\la)\in R(A)\MOF$. Therefore, (see \eg 
\cite[Theorem~3]{Frisktilting} and \cite[Theorem~1.6]{AHLU}), 
$\mathtt{R}Q$ has a proper costandard flag, hence it is tilting.     
\end{proof}

For any abelian category $\cC$ we denote by $\cD^b(\cC)$ its bounded
derived category. If $\cC=A\MOF$ we set $\cD^b(A)=\cD^b(\cC)$.
If the opposite is not explicitly stated, by a ``functor'' we always
mean a covariant functor. We use the standard notation like $\cL F$, $\cR G$,  
$\cL_i F$, $\cR^i G$ etc.\ to denote left derived and right derived functors 
and their $i$-th cohomology functors.

For the sake of completeness we mention the following fact

\begin{prop}
\label{Ringel2}
Let $A$ and $B$ be standardly stratified such that tilting modules
are also cotilting. Let $F:A\MOF\rightarrow B\MOF$ be a (covariant) right-exact 
functor with right adjoint $G$. Assume that $F$ defines an equivalence 
\begin{displaymath}
\cF(\overline{\Delta}^A)\cong\cF(\overline{\nabla}^B).
\end{displaymath}
Then the following hold
\begin{enumerate}[(1)]
\item $F$ maps projective modules to tilting modules and tilting modules to
injective modules. In fact, $F$ defines equivalences (with inverse $G$) of 
the corresponding additive subcategories.
\item If $A$ has finite global dimension then $B$ has 
finite global dimension as well, moreover, 
$\cL F:\cD^b(A)\rightarrow\cD^b(B)$ is an equivalence with
inverse $\cR G$. 
\item $B$ is the Ringel dual of $A$.   
\end{enumerate}
\end{prop}

\begin{proof}
Let $P\in A\MOF$ be projective, then $FP\in\cF(\overline{\nabla}^B)$ by assumption and
\begin{displaymath}
\EXT^1_B(FP,\overline{\nabla}(\la))\cong\EXT^1_B(FP,FF^{-1}\overline{\nabla}(\la))
\cong\EXT^1_A(P,F^{-1}\overline{\nabla}(\la))=0
\end{displaymath}
for any proper costandard module $\overline{\nabla}(\la)$. 
Hence $FP\in\cF(\Delta^B)$ and is therefore
tilting. If $X$ is tilting, hence cotilting, then $X\in\cF(\overline{\Delta}^A)$.
Therefore $FX\in\cF(\overline{\nabla}^B)$ and
\begin{eqnarray}\label{incsplits}
\EXT^1_B(\overline{\nabla}^B(\la),FX)&\cong&
\EXT^1_B(FF^{-1}\overline{\nabla}^B(\la), FX) \nonumber\\
&\cong&\EXT^1_A(F^{-1}\overline{\nabla}^B(\la),X)=0
\end{eqnarray}
for any proper costandard module $\overline{\nabla}^B(\la)$,
since $F^{-1}\overline{\nabla}^B(\la)\in\cF(\overline{\Delta}^A)$ and
$X\in \cF(\nabla^A)$. If we now choose an inclusion of $FX\in 
\cF(\overline{\nabla}^B)$ in its injective hull, then the cokernel is 
contained in $\cF(\overline{\nabla}^B)$ and the inclusion splits 
because of \eqref{incsplits}. This means that $F X$ is injective and
the first part follows. We have $\cR G\,\cL F\cong\ID$ on 
projectives and $\cL F\,\cR G\cong\ID$ on injectives. This implies that
the global dimension of $B$ is finite and then the second statement follows. 
The fact that $B$ is the Ringel dual of $A$ is then clear from the definitions. 
\end{proof}

\subsection{Double centraliser property}
We claim that, given a standardly stratified algebra $A$, there is 
{\it always} some tilting module $X$ such that we have a double centraliser
property,  $A\cong\END_{\END_A(X)}(X)$. This relies on the following 

\begin{prop}
\label{whichtilting}
Let $A$ be standardly stratified and let $R=R(A)$ be its Ringel dual. Let $P$ be
the projective cover of the characteristic tilting module $T$ in $\op{mod-}R$. Then
there is an exact sequence $0\rightarrow A\rightarrow Q\rightarrow \COKER\rightarrow 0$,
where $Q=\mathtt{R}^{-1}P$ (see Proposition~\ref{Ringeldual}) is tilting and 
$\COKER\in\cF(\Delta^A)$.  
\end{prop}

\begin{proof}
Since $P,T\in\cF(\Delta^R)$, the kernel $K$ of the surjection between $P$ and
$T$ is contained in $\cF(\Delta^R)$ (\cite[Theorem 1.6 (i)]{AHLU} and
\cite[Theorem~3]{Frisktilting}). Applying the inverse of 
the Ringel duality functor (which is defined on $\cF(\Delta^R)$) 
we get the short exact sequence 
\begin{displaymath}
0\rightarrow
A\rightarrow Q\rightarrow \COKER\rightarrow 0,
\end{displaymath}
where $\COKER\in\cF(\Delta)$ by Proposition~\ref{Ringeldual}.
\end{proof}

\begin{cor}
\label{doublecentralizer}
There exists a (basic) tilting module $X$ such that we have an isomorphism,
$A\cong\END_{\END_A(X)}(X)$. 
\end{cor}

\begin{proof}
Let $Y$ be a tilting module such that we have an inclusion $\COKER\hookrightarrow
Y$ (the existence follows from \cite[Theorem~5.4]{AR1}). 
Put $X:=(Q\oplus Y)_{basic}$, then there exists an exact sequence, 
$0\rightarrow \COKER\rightarrow X^n$,
satisfying the assumptions of \cite[Theorem 2.8]{KSX}. Hence the double
centraliser property follows.
\end{proof}

\begin{remark}
{\rm
One can show that there exists a {\em minimal} basic tilting module $Y$ with the 
following property: any $M\in \cF(\Delta)$ embeds into $Y^m$ for some $m$.
Here {\em minimal} means that every other tilting module with the latter property
has $Y$ as a direct summand. However, it is not clear whether there exists a
minimal basic tilting module $Y$, with respect to which
one has the double centraliser property. It is the case in all the
examples we know, in particular in  the Examples~\ref{examples}.

In general, it could happen that $X$ is already the characteristic tilting module,
and the statement of Corollary~\ref{doublecentralizer} 
is not very useful. As an example we refer to
\cite[Example $A_1$]{KK} where the algebra $A$ is given by all $3\times 3$ 
upper triangular matrices over some field $\Bbbk$ with the matrix idempotents 
$e_1$, $e_2$, $e_3$ and the quasi-hereditary structure given by the ordering 
$1<2<3$. The same algebra, but with the quasi-hereditary structure given by the
reversed order (see \cite[Example $A_2$]{KK}) provides also an example, where 
$X$ is not contained in $\op{Add}(Q)$, the additive category generated by $Q$. 
In particular, we do not have the double centralizer property with respect to $Q$.
}
\end{remark}

\subsection{Double centralizer and projective-injective modules}
On the other hand, under the assumptions and notation of 
Proposition~\ref{whichtilting} we have the following nice situation, where
projective-injective modules play a crucial role.

\begin{cor}\label{nice}
If the injective hull of any standard module is contained in
$\op{Add}(Q)$, then the following holds:
\begin{displaymath}
A\cong\END_{\END_A(Q)}(Q)\cong\END_{\END_A(Q_{basic})}(Q_{basic}).
\end{displaymath}
\end{cor}

\begin{proof}
If the injective envelope of any standard module is contained in $\op{Add}{Q}$,
the assumptions of \cite[Theorem~2.8;~Theorem~2.10]{KSX} are satisfied and the 
statement follows. 
\end{proof}

As interesting examples we have the following:

\begin{exs}\label{examples}
{\rm
In the examples which follow we illustrate the use of Proposition~\ref{whichtilting} 
and Corollary~\ref{doublecentralizer}, in particular, we
explicitly describe the modules $Q$ and $X$ which appear in the double centraliser
statements.
\begin{enumerate}
\item
\label{excatO}
Let $A$ be such that $A\MOF$ is equivalent to an integral block of
the Bernstein-Gelfand-Gelfand category $\cO$ for some semi-simple complex Lie
algebra $\mg$ (see \cite{BGG}). The algebra $A$ is equipped with the usual 
quasi-hereditary structure (given by the Bruhat order and the Verma modules as
standard modules). 
In this case we have exactly one indecomposable projective-injective module, namely
the projective cover $P(w_0)$  of the unique simple standard (or Verma) module
in this block. Moreover, $A$ is Ringel self-dual (\cite[Theorem 5.12 and
Bemerkung 2.4 (3)]{SoKac}). The projective cover of a tilting 
module is a direct sum of $P(w_0)$'s. Via Ringel duality we get an inclusion 
\begin{equation}
\label{step1}
  i:A\hookrightarrow Q,
\end{equation}
where $Q=P(w_0)^n$ for some positive integer $n$. The cokernel of this inclusion has (by
Proposition~\ref{whichtilting}) a standard (or Verma-) flag. Hence there is an
exact sequence of the form 
\begin{equation}
\label{step2}
  0\rightarrow A\longrightarrow Q\longrightarrow Q^m,
\end{equation}
for some positive integer $m$. We could take  $X=Q_{basic}=P(w_0)$ and get the famous double centraliser theorem of 
Soergel (\cite{Sperv}, see also \cite[Theorem 3.2]{KSX}), namely
$A\cong\END_{\END_A(X)}(X)$.  
\item 
\label{exparabO}
Let $A^\p$ be such that $A^\p\MOF$ is equivalent to an integral block
of some parabolic category $\cO^\p$ in the sense of \cite{RC} (see also
Section~\ref{parabolic}) with the usual quasi-hereditary structure. 
Then $A^\p$ is Ringel-self-dual (see \cite{SoKac} or Proposition~\ref{Serreparab} 
below). The self-dual  projective modules are exactly the
summands occurring in the injective hulls of standard modules
(\cite{Irself}), they are also exactly the summands occurring in the
projective cover of tilting modules. This means, we have an embedding of the
form~\eqref{step1} and then an exact sequence of the form~\eqref{step2}, where
$Q$ is a direct sum of projective-injective modules. If we set $X=Q_{basic}$ the sum over 
(a system of representatives for the isomorphism classes
of) all indecomposable projective-injective modules we get the double
centraliser property $A^\p\cong\END_{\END_{A^\p}(X)}(X)$ (this is proved in
\cite[Theorem 10.1]{Stquiv}).  
\item 
\label{exHC}
Let $\mg$ be a semisimple complex Lie algebra. Let $\cH$ be the category
of Harish-Chandra bimodules for $\mg$, that is the category of
$\mg$-bimodules which are of finite length and locally finite for the
adjoint action of $\mg$ (see for example~\cite{BG} or \cite[Section 6]{Ja2}). 
The category $\cH$ decomposes into blocks
$_\la\cH_\mu$. A bimodule $X\in\cH$ is contained in the block $_\la\cH_\mu$ if it
is annihilated by $(\ker\chi_\la)^n$ from the left and by $(\ker\chi_\mu)^n$ from
the right for some positive integer $n$, where $\ker\chi_\la$ is the annihilator of the
Verma module with highest weight $\la$. The category $_\la\cH_\mu$ does not have projective
objects, however, we get enough projectives (see \eg \cite[6.14]{Ja2}) if we consider 
the full subcategory $^{}_\la\cH_\mu^1$ of
$^{}_\la\cH_\mu$ given by all bimodules which are annihilated by
$\KER\chi_\mu$ from the right hand side.     
Let $A_\la^\mu$ be such that $A_\la^\mu\MOF\cong{} ^{}_\la\cH_\mu^1$, where $\la$ 
and $\mu$ are integral. Then $A_\la^\mu$ is standardly stratified (it is not
quasi-hereditary in general) and contains a unique
indecomposable projective-injective module (see
\cite[Corollary~2]{KoMa}). Later (Proposition~\ref{SerreHC}) we give a new proof
for the fact that $A_\la^\mu$ is Ringel self-dual (see \cite[Theorem~3]{FKM1}
for the original argument). As in category
$\cO$, the projective cover of a tilting module is projective-injective, and hence
$Q$ becomes a direct sum of copies of the unique self-dual indecomposable projective 
module. The injective hulls of standard
modules are projective as well. Hence we could take $X=Q_{basic}$, 
the indecomposable projective-injective module and get the double centraliser $A_\la^\mu\cong\END_{\END_{A_\la^\mu}(X)}(X)$.   
\item 
\label{exCherednik}
Let $A$ be such that $A\MOF\cong\cO(H_c)$, the category $\cO$ for
some rational Cherednik algebra $H_c=H_{0,c}$ as considered for example in
\cite{Guay} or \cite{GGOR}.  The projective-injective modules are exactly the
summands occurring in the injective hulls of standard modules
(\cite[Proposition 5.21]{GGOR}), they are also exactly the summands occurring in the
projective covers of tilting modules. Hence, $Q$ is a direct sum of projective-injectives 
and then we could take $X=Q_{basic}$ to be the sum over all indecomposable projective 
tilting modules. This is the double centraliser property from \cite[Theorem 5.16]{GGOR}.
\item \label{exquiver}
Quite often there are double centraliser properties with respect to tilting modules, which
do not have to be projective or injective. In the following examples the tilting module 
$X$ is neither projective nor injective: 
Let $\mathfrak{Q}$ be a finite quiver with vertices $\{1,\dots,n\}$. Assume it is
directed, that is an arrow from $i$ to $j$ exists only if $i>j$. Let $A=A(\mathfrak{Q})$ 
be the corresponding path algebra and $D$ be its dual extension, that is the algebra
$A\otimes_{\Bbbk}A^{opp}$ with the relations $(\mathrm{rad}A^{opp})(\mathrm{rad}A)=0$
(see \eg \cite{DX}). Then $D$ is quasi-hereditary  with respect to the natural
order on $\{1,\dots,n\}$. One can show that  there is a double centraliser
property with respect to the tilting module $X=Q_{basic}=\oplus_{i}T(i)$, where the sum runs
over all sources of $\mathfrak{Q}$. It is also easy to see that $X$ is neither
injective nor projective in general. 
\end{enumerate}
}
\end{exs}
\begin{remark}
{\rm Let $A$ and $X$ be as in the examples above, then we could define 
\begin{eqnarray*}
  \mV: A\MOF&\longrightarrow& \END_A(X)\MOF\\
  M&\longrightarrow&\HOM_A(X,M).
\end{eqnarray*}
The double centraliser property can be reformulated as: The functor $\mV$ is fully faithful on projective modules, \ie $\mV$ induces an
  isomorphism
  \begin{eqnarray*}
    \HOM_A(P_1,P_2)&\cong&\HOM_{\END_A(X)}(\mV P_1, \mV P_2)
  \end{eqnarray*}
for all projective modules $P_1$ and $P_2$. 

Another easy consequence from the definitions is the following: 
The functor $\mV$ is fully faithful on tilting modules, \ie $\mV$ induces an
  isomorphism
  \begin{eqnarray}
\label{faithtilt}
    \HOM_A(T_1,T_2)&\cong&\HOM_{\End_A({X})}(\mV T_1, \mV T_2)
  \end{eqnarray}
for all tilting modules $T_1$ and $T_2$. 

\begin{proof}
If $\Hom_A(T_1,K)=0=\Hom_A(K,T_2)$ for any $K\in A\MOF$ such that $\mV K=0$ then 
  \begin{eqnarray*}
    \HOM_A(T_1,T_2)&\cong&\HOM_{\End_A(X)}(\mV T_1, \mV T_2)
  \end{eqnarray*}
since $\mV$ is a quotient functor (see \cite{Gabriel}). 
All the composition factors in $K$ are annihilated by $\mV$. On the other
hand, none of the composition factors in the head of $T_1$ and in the socle of
$T_2$ is annihilated by $\mV$. This proves the statement. 
\end{proof}
}
\end{remark}
   
\section{Serre functors}\label{section2}

The aim of the present section is to develop an effective machinery to 
describe  Serre functors for the categories appearing in the examples above, where 
the algebra is not
given explicitly. Let $\cC$ be a $\Bbbk$-linear additive category with finite 
dimensional homomorphism spaces. A {\it right Serre functor} is an additive 
endofunctor $F$ of $\cC$ together with isomorphisms 
\begin{eqnarray}
\label{Serreright}
  \Psi_{X,Y}: \HOM_\cC(X,F Y)\cong\HOM_\cC(Y,X)^\ast,
\end{eqnarray}
 natural in $X$ and $Y$. Here, $\ast$ denotes the ordinary duality for vector spaces.
Right Serre functors satisfy the following properties:
\begin{itemize}
\item Two right Serre functors are isomorphic 
(see \cite[Lemma~I.1.3]{RvB}).
\item If $\epsilon$ is an auto-equivalence of $\cC$ and $F$ is a right Serre functor 
then $\epsilon F\cong F\epsilon$. (It follows directly from the definitions that 
$\epsilon F\epsilon^{-1}$ is a right Serre functor, hence it must be isomorphic to $F$).    
\end{itemize}
A right Serre functor is a {\it Serre functor} if it is an auto-equivalence of
$\cC$. By general results (see \cite{BoKa}), for any finite dimensional algebra $A$ of
finite global dimension, there is a Serre functor $\mS$ for the bounded
derived category $\cD^b(A)$, more precisely $\mS\cong\cL H$, where
$H=A^*\otimes_A\bullet$ (\cite[Example~3.2(3)]{BoKa}). In the literature, 
the functor $H$ is often called the
{\it Nakayama functor} (see e.g. \cite[page 37]{Happel}). This is because $H\cong
\HOM_A(\bullet, A)^\ast$. 

Recall that for any abelian category $\cC$ we denote by $\cD^b(\cC)$ its bounded
derived category. If $\cC=A\MOF$ we set $\cD^b(A)=\cD^b(\cC)$. We use the standard notation like $\cL F$, $\cR G$,  
$\cL_i F$, $\cR^i G$ etc.\ to denote left derived and right derived functors 
and their $i$-th cohomology functors. Let also $\mathcal{D}_{\mathrm{perf}}(A)$
denote the full subcategory of $\mathcal{D}^b(A)$, consisting of {\em perfect} 
complexes (\ie of those complexes which are quasi-isomorphic to bounded complexes 
of projective $A$-modules).

In order to be able to describe more explicitly the Serre functors for some of the examples 
mentioned above we will need effective tools to detect Serre
functors. Recall that a finite-dimensional algebra, $A$, is called {\em self-injective}
provided that $A\cong A^*$ as left $A$-modules; and {\em symmetric} provided that
$A\cong A^*$ as $A$-bimodules. We start with the following easy observation

\begin{lemma}\label{lperfect}
Let $A$ be a finite-dimensional self-injective algebra. Then $\mathcal{L}H$ is a Serre functor of $\mathcal{D}_{\mathrm{perf}}(A)$, 
moreover, $\mathcal{L}H\cong\ID$ if and only if $A$ is symmetric.
\end{lemma}

\begin{proof}
Let $\mathcal{P}^{\bullet}$ be a bounded complex of projective $A$-modules. Then
we have that $\mathcal{L}H\mathcal{P}^{\bullet}=H\mathcal{P}^{\bullet}$ is a bounded complex 
of injective $A$-modules by the definition of $H$. Since $A$ is self-injective
we have $H\mathcal{P}^{\bullet}\in\mathcal{D}_{\mathrm{perf}}(A)$. That in this 
case $\mathcal{L}H$ is a Serre functor is proved for example in  
\cite[Proposition~20.5.5(i)]{Ginzburg}. Finally, the last statement follows from
the definition of a symmetric algebra.
\end{proof}

\begin{defn}
{\rm 
Given an algebra $A$ and a projective-injective module $Q$, we call
$Q$ {\it good} if the socle of $Q$ is isomorphic to the head of
$Q$. (Equivalently, if $Q\cong\oplus_{\la\in\Lambda'} P(\la)$ for some
$\Lambda'\subset\Lambda$ then $Q\cong \oplus_{\la\in\Lambda'} I(\la)$.)} 
\end{defn}

If $A$ has a duality which preserves simple modules, any projective-injective
module is automatically good.    

\begin{remark}\label{rem504}
{\rm
In the following we will also use double centraliser properties for the opposed
algebra $A^{opp}$. Let $I$ be a basic injective $A$-module. It is easy to see
that the existence of an exact sequence of the form 
\begin{eqnarray}
\label{ess}
  Q_2\rightarrow Q_1\rightarrow I\rightarrow 0 
\end{eqnarray}
for some projective-injective $A$-modules $Q_1$, $Q_2$, is equivalent to the
requirement that $A^{opp}$ has a double centraliser property with respect to a
projective-injective module. Indeed, the double
centraliser property for $A^{opp}$ is equivalent to the existence of an exact
sequence of the form 
\begin{equation}
0\rightarrow A^{opp}\longrightarrow X_1'\longrightarrow  X_2'
\end{equation}
for some projective-injective modules $X_1'$, $X_2'$. Applying the usual
duality $\HOM_{\Bbbk}(\bullet, \Bbbk)$ we get an exact sequence
\begin{equation}\label{pereqex1}
X_2\to X_1\to I\to 0,
\end{equation}
where $I$ is the injective cogenerator of $A\mathrm{-mod}$ and $X_1$, $X_2$
are projective-injective.
}
\end{remark}

\subsection{A characterisation of Serre functors}
The following result provides a tool to detect Serre functors:

\begin{theorem}\label{nnnnnnn}
Let $A$ be a finite dimensional $\Bbbk$-algebra of finite global dimension. 
Assume that a basic projective-injective $A$-module is good and both, $A$ and
$A^{opp}$, have the double centraliser property with respect
to a projective-injective module. Let $F:A\MOF\rightarrow A\MOF$ be a right exact functor. Then $\cL F$ is a 
Serre functor of $\mathcal{D}^b(A)$ if and only if the following conditions 
are satisfied:
\begin{enumerate}[(a)]
\item\label{nnnnnnn.1}
Its left derived functor $\cL F:\cD^b(A)\rightarrow \cD^b(A)$ is an auto-equivalence.
\item\label{nnnnnnn.2}
$F$ maps projective $A$-modules to injective $A$-modules.
\item\label{nnnnnnn.3} 
$F$ preserves the full subcategory $\mathcal{PI}$ of $A\mathrm{-mod}$, consisting of
all pro\-jec\-tive-\-in\-jec\-tive modules, and the restrictions of $F$ and $H$ to
$\mathcal{PI}$ are isomorphic. 
\end{enumerate}
\end{theorem}

\begin{proof}
Let $Q$ be a good basic projective-injective $A$-module. We know that $\cD^b(A)$ has a Serre functor, $\mS$, and $\mS\cong\cL H$,
where $H=A^*\otimes_A\bullet$. By definition, $H$ satisfies \eqref{nnnnnnn.1}
and \eqref{nnnnnnn.2} and preserves $\mathcal{PI}$, because $Q$ is good. Hence $H$ satisfies \eqref{nnnnnnn.3}.  

Now let $F:A\MOF\rightarrow A\MOF$ be a right exact functor, satisfying
\eqref{nnnnnnn.1}--\eqref{nnnnnnn.3}. We claim that $F$ and $H$ are isomorphic
when restricted to the category of injective $A$-modules. Indeed, the double
centraliser property for $A^{opp}$ gives us an exact sequence,
\begin{equation}\label{pereqex}
X_2\to X_1\to I\to 0,
\end{equation}
where $I$ is the injective cogenerator of $A\mathrm{-mod}$ and $X_1,X_2\in
\mathcal{PI}$ (see Remark~\ref{rem504}).
Let $\psi:F\to H$ be the isomorphism, given  by \eqref{nnnnnnn.3}. Applying
$F$ and $H$ to \eqref{pereqex} and using \eqref{nnnnnnn.3} we obtain the following
diagram with exact rows, where the square on the left hand side commutes,
inducing an isomorphism, $\psi_I$, as indicated:
\begin{displaymath}
\xymatrix{
F(X_2)\ar[rr]\ar[d]_{\psi_{X_2}} & & F(X_1)\ar[rr]\ar[d]_{\psi_{X_1}} & & 
F(I)\ar[rr]\ar@{.>}[d]_{\psi_I} & & 0\\ 
H(X_2)\ar[rr] & & H(X_1)\ar[rr] & & H(I)\ar[rr] & & 0\\ 
}
\end{displaymath}
By standard arguments, it defines an isomorphism of functors, $F\cong H$, when
restricted to the full additive category of injective $A$-modules. Since $\cL
F$ is an auto-equi\-va\-lence, we have $\cL F\,\mS\cong\mS\,\cL F$. 
As projectives are acyclic for right exact functors, we get an isomorphism, $\cL
F \, H\cong\mS \, F$, when restricted to the full additive subcategory given
by projectives. Taking the $0$-th
homology we get an isomorphism of functors 
\begin{equation}
  \label{eq:FHHF}
F\,H\cong H\, F  
\end{equation}
when restricted
to the full additive category of projective $A$-modules. Since the functors $F$ and
$H$ are right exact, we only have to deduce that $F\cong H$ on the category of
projectives. We already know that $F$ and $H$ are invertible on
$\mathcal{PI}$, hence we can fix isomorphisms
$\alpha:\END_A(Q)\cong\END_A(FQ)$ and $\beta:\END_A(Q)\cong\END_A(HQ)$. When
restricted to $\mathcal{PI}$, we have $F\cong\ID^\alpha$ and
$H\cong\ID^\beta$, where $\ID^\alpha$ and $\ID^\beta$ denote the identity
functors, but with the $\END_A(Q)$-action twisted by $\alpha$ or $\beta$
respectively. Since both, $F$ and $H$, are right exact, they uniquely extend to
functors on $\op{mod}$-$\END_A(Q)$, the latter being realized as the full subcategory 
$\cC$ of $A\MOD$ given by all modules, having a presentation of the
form~\eqref{pereqex1} (see \cite[Section~5]{Au}). From the explicit description above, we obtain
that both $F$ and $H$ are invertible as endofunctors of $\cC$. As both, $H$ and $F$, 
map projectives to injectives and $F\cong H$ on injectives we get, together
with ~\eqref{eq:FHHF}, isomorphisms of functors $F^2\cong
H\,F\cong F\, H$ when restricted to the full additive category of projective
$A$-modules. This gives then rise to an isomorphism, $F\cong H$, since $F$ is invertible on $\cC$. So, we are done.
\end{proof}

\begin{prop}\label{cperfect}
Let  $A$ be a finite dimensional $\Bbbk$-algebra of finite global
dimension. Assume there is a good basic projective-injective module $Q$ and set
$B=\mathrm{End}_{A}(Q)$. Then the algebra $B$ is symmetric if and
only if the restriction of the Serre functor for $\cD^b(A)$ to $\mathcal{PI}$ is
the identity functor. 
\end{prop}

\begin{proof}
Let $\mS$ be the Serre functor for $\cD^b(A)$. Then $\mS$ obviously preserves
$\mathcal{PI}$, because $Q$ is good, and hence it also preserves the (homotopy) category of bounded 
complexes of projective-injective $A$-modules. Moreover, it induces a Serre functor
on this category. By \cite[Section~5]{Au}, the latter one is equivalent to the
category $\mathcal{D}_{\mathrm{perf}}(B)$. The statement now follows from
Lemma~\ref{lperfect}.
\end{proof}

\subsection{Serre functors via partial coapproximation}
In this subsection we want to show that double centraliser properties with
respect to projective-injective modules quite often make it possible to
describe the Serre functor in terms of partial coapproximations. 

For the remaining section let $A$ be a finite dimensional algebra of finite
global dimension. 
Let $Q\in A\MOD$ be a projective module. For any module $M$ let $M_Q$ be the
trace of $Q$ in $M$ (ie. $M_Q$ is the smallest submodule of $M$ such that
$\HOM_A(Q,M/M_Q)=0$). Dually let $M^Q$ be the smallest quotient of $M$ such
that $\HOM_A(Q,M)=\HOM_A(Q,M^Q)$. 

Associated with $Q$, there is a right exact functor $\CA_Q:A\MOD\rightarrow A\MOD$ called the
{\it partial coapproximation with respect to $Q$} (for details we refer for
example to \cite[2.5]{KM}). It sends a projective
module $P$ to $P_Q$. Note that if $f:P\rightarrow P'$ is a morphism between 
projective modules, then it induces a morphism, $\CA_Q(f):\CA P\rightarrow\CA(P')$. 
These assignments can be extended uniquely to a right exact endofunctor $\CA_Q$ 
of $A\MOD$. For an arbitrary module $M\in A\MOD$, the module $\CA_Q
M$ can be constructed in the following way: We choose a short exact sequence
$K\hookrightarrow P\twoheadrightarrow M$, where $P$ is projective. Then
\begin{displaymath}
\CA_Q M\cong \left(P/K_Q\right)_Q,
\end{displaymath}
in other words $\CA_Q M$ is obtained from $M$ by first maximally extending $M$
using simple modules, which do not occur in the top of $Q$, and afterwards deleting
all occurrences of such modules in the top part.

\begin{lemma}
\label{essential}
Let  $A$ be a finite dimensional $\Bbbk$-algebra. Assume, $A^{opp}$ has the
double centraliser property with respect to a projective-injective module. Let
$Q$ be a basic projective-injective $A$-module. Let $\lambda\in\Lambda$. Then the following holds: If
$P(\lambda)_Q\cong I(\lambda)^Q$ then $(\CA_Q)^2(P(\lambda))\cong I(\lambda)$.
\end{lemma}

\begin{proof}[Proof of Lemma~\ref{essential}]
  We have 
  \begin{eqnarray*}
    (\CA_Q)^2(P(\lambda))&\cong&\CA_Q (P(\lambda)_Q)\\
               &\cong&\CA_Q (I(\lambda)^Q)\\
               &\cong&I(\lambda).
  \end{eqnarray*}
Here, only the last isomorphism needs some explanation. If $P$ is the
projective cover of $I(\lambda)^Q$ then the natural surjection from $I(\lambda)$ onto
$I(\lambda)^Q$ lifts to a map, $f:P\rightarrow I(\lambda)$. From the definition of 
$I(\lambda)^Q$ and \eqref{ess} it follows that $f$ is surjective.
The double centraliser property for $A^{opp}$ (see Remark~\ref{rem504}) also implies that any composition
factor in the head of the kernel of $f$ is not annihilated by
$\HOM_A(Q,\bullet)$. Hence the desired isomorphism follows.  
\end{proof}

The following theorem describes a situation, where the double centraliser
property with respect to a basic projective-injective module $Q$, the
description of the Serre functor via partial coapproximation, and the symmetry
of the endomorphism ring of $Q$ are nicely connected. Later on we will see
that this setup applies to
all the different versions of category $\cO$ mentioned in the Examples~\ref{examples}. 

\begin{theorem}
\label{Serrecoapprox}
Let  $A$ be a finite dimensional $\Bbbk$-algebra of finite global
dimension. Let $Q$ be a basic projective-injective $A$-module. Assume, $Q$ is
good and both, $A$ and $A^{opp}$ have the double centraliser property
with respect to some projective-injective module. Consider the functors
$\mV=\HOM_A(Q,\bullet):A\MOD\rightarrow\op{mod}$-$\END_A(Q)$ and
$H=A^\ast\otimes_A\bullet$. Then the following assertions are equivalent:
\begin{enumerate}[(i)]
\item \label{i} $\mV\cong\mV H$,
\item \label{ii} $H\cong(\CA_Q)^2$,
\item \label{iii} $\END_A(Q)$ is symmetric.
\end{enumerate}
In either of these cases, the Serre functor for $\cD(A)^b$ is $\cL((\CA_Q)^2)$.
\end{theorem}

\begin{proof} 
Obviously, if~\eqref{ii} holds then $\cL((\CA_Q)^2)$ is the Serre functor for
$\cD^b(A)$. It is left to show that the three cases are equivalent.\\
{$\eqref{i}\Rightarrow\eqref{ii}:$} Let us assume $\mV\cong\mV H$. Let $P$ be
a projective module. By the assumed double centraliser property for $A$
and $A^{opp}$ (see Remark~\ref{rem504}) we have natural isomorphisms 
\begin{displaymath}
\begin{array}{rcl}
\HOM_A(P,P)&\cong&\HOM_{\END_A(Q)}(\mV P, \mV P)\\
&\cong&\HOM_{\END_A(Q)}(\mV H P, \mV P)\\
&\cong&\HOM_A(H P, P).
\end{array}
\end{displaymath}
(for the last 
isomorphism we refer to the proof of \eqref{faithtilt}). The identity map
in $\END_A(P)$ gives rise to a natural morphism, $HP\rightarrow P$, identifying
$(HP)^Q$ and $P_Q$. Since $H$ maps the projective cover of any simple module to its
injective hull, we are in the situation of Lemma~\ref{essential}. In
particular, $(\CA_Q)^2$ sends an indecomposable projective module to the corresponding
indecomposable injective module. Let $G$ be the right adjoint functor to
$\CA_Q$ (this is the functor of partial approximation with respect to $Q$, see
\cite[2.5]{KM}). We
have the adjunction morphism $\ID\rightarrow G^2(\CA_Q)^2$ which we know is an
isomorphism on projective-injective modules. From the double centraliser
property we get that this adjunction morphism is injective on all projective
modules. Since $(\CA_Q)^2 P$ is isomorphic to the corresponding injective module, we have
$G^2(\CA_Q)^2\cong P$. In particular, $G^2(\CA_Q)^2\cong\ID$ when restricted to the
additive subcategory given by projective modules. Dually,
$(\CA_Q)^2G^2\cong\ID$ when restricted to the additive category given by
injective modules. Since $A$ has finite global dimension, $\cL((\CA_Q)^2)$ defines
an auto-equivalence of the derived category $\cD^b(A)$ with inverse $\cR
G^2$. From our assumption we have $\mV(\CA_Q)^2\cong\mV\cong\mV H$. Therefore, $(\CA_Q)^2\cong H$ on the
additive subcategory given by all projective-injective modules. Hence
$(\CA_Q)^2$ satisfies the assumptions of Theorem~\ref{nnnnnnn}. It follows in particular,
$H\cong (\CA_Q)^2$.
\\
{$\eqref{ii}\Rightarrow\eqref{iii}:$} The definition of $\CA_Q$ implies that it induces the identity functor on
the category of projective-injective $A$-modules. Hence $\END_A(Q)$ is symmetric by Proposition~\ref{cperfect}.\\
{$\eqref{iii}\Rightarrow\eqref{i}:$} We assume that $B=\END_A(Q)$ is
symmetric. From Lemma~\ref{lperfect} we have that the Serre functor of
$\cD_{perf}(B)$ is isomorphic to the identity functor. On the other hand, the Serre functor of
$\cD^b(A)$ induces a Serre functor on the category of bounded complexes of
projective-injective $A$-modules. (Note that this category is preserved by the
Serre functor, since $Q$ was assumed to be good.) Altogether, when 
restricted to the category of projective-injective modules, the functor $H$ is isomorphic to the
identity functor. This provides the following sequence of natural isomorphisms for
any projective $A$-module $P$: 
\begin{eqnarray*}
  \mV H P&\cong&\HOM_A(Q, HP)\\
&\cong&\HOM_A(HQ, HP)\\
&\cong&\HOM_A(Q,P)\\
&\cong&\mV P.
\end{eqnarray*}
(For the penultimate isomorphism we used that $H$ defines an auto-equivalence of
$\cD^b(A)$, hence it is in particular fully faithful on projectives.) Thus we get an
isomorphism of functors $\mV H\cong \mV$ when restricted to the category of
projective modules. Since the involved
  functors are right exact, the isomorphism extends to an isomorphism of
  functors  $\mV H\rightarrow\mV$.
\end{proof}

\section{Applications}\label{secappl}

\subsection{Bernstein-Gelfand-Gelfand category $\cO$}\label{secappl.1}
Let $\mg$ be a semisimple complex Lie algebra with a fixed Borel subalgebra
$\mathfrak{b}$ containing the fixed Cartan subalgebra $\mh$. Let $\cO$ be the
corresponding BGG-category (see \cite{BGG}). Let $W$ denote the Weyl group of
$\mg$ with longest element $w_0$. For any weight $\la\in\mh^\ast$ let $W_\la$
be the stabiliser $W_\la=\{w\in W\mid w\cdot\la=\la\}$, where
$w\cdot\la=w(\la+\rho)-\rho$ and $\rho$ is the half-sum of positive roots.   
For $\mu\in \mh^\ast$ let $\Delta(\mu)$ be the
Verma module with highest weight $\mu$. For $\la\in\mh^\ast$, a dominant
and integral weight, we consider the block $\cO_\la$, containing the Verma modules 
$\Delta(\mu)$, where $\mu\in W\cdot \la$. Let $L(\mu)$ be the simple quotient of 
$\Delta(\mu)$ and $P(\mu)$ its projective cover. For any $w\in W$, there is a 
{\it twisting functor} $T_w:\cO\to\cO$ (given by tensoring with some ``semi-regular 
bimodule''), see \cite{AL}, \cite{KM} or \cite{AS} for a precise definition. 
Let $\op{d}$ be the duality on $\cO$. We denote by $G_w$ the 
right adjoint functor of $T_w$. We have  $G_w\cong\op{d}T_w\op{d}$ 
(see \cite[Section~4]{AS}). 

If $\la$ is regular, and $s$ is a simple reflection, we denote by $C_s$ Irving's {\it shuffling functor} defined as
taking the cokernel of the adjunction morphism between the identity functor
and the translation $\theta_s$ ``through the $s$-wall'' (\cite[Section~3]{GJ}). Let
$w_0=s_{i_1}s_{i_2}\cdots s_{i_r}$ be a reduced expression, then we define
$C_{w_0}= C_{s_{i_r}}C_{s_{i_{r-1}}}\cdots C_{s_{i_1}}$. Up to isomorphism,
this does not depend on the chosen reduced expression (see \eg \cite[Lemma~5.10]{MS}).

\begin{prop}\label{prcato}
Let $A=A_\la$ such that $A_\la\MOF\cong\cO_\la$ for some integral block
$\cO_\la$. 
\begin{enumerate}[(1)]
\item The functor $\cL\,(T_{w_0})^2:\cD^b(\cO_\la)\rightarrow\cD^b(\cO_\la)$
  is a Serre functor.    
\item If $\la$ is regular, then
  $\cL\,(C_{w_0})^2:\cD^b(\cO_\la)\rightarrow\cD^b(\cO_\la)$ is a Serre functor. In particular,
  $\cL\,(T_{w_0})^2\cong\cL\,(C_{w_0})^2$. 
\end{enumerate}
\end{prop}

\begin{proof}
We want to verify the assumptions of Theorem~\ref{nnnnnnn} for $A=A_\la$ and
$F= T_{w_0}^2$ considered as an endofunctor of $A\MOD$.
  
Because of the existence of a duality on $A$ we have $A\cong A^{opp}$ and,
as we have already mentioned in the introduction, $A$ has a double centraliser
property with respect to the good basic projective-injective module
$P(w_0\cdot\la)$ (see \cite[Struktursatz]{Sperv}). If $\la$ is regular, the endomorphism algebra of
the latter is the coinvariant algebra associated
with $W$. If $\la$ is singular then this endomorphism
ring is isomorphic to the subalgebra of $W_\la$-invariants in the coinvariant algebra
(\cite[Endomorphismensatz]{Sperv}). In any case, the resulting algebra is
symmetric. 
Consider now $(T_{w_0})^2:\cO_\la\rightarrow\cO_\la$. This functor is both right exact and 
additive by definition. It's derived functor defines a self-equivalence of
$\mathcal{D}^b(\cO_\la)$ by \cite[Corollary~4.2]{AS} for the regular case; the
singular case follows by translation, since twisting functors commute
naturally with translation functors (\cite[Theorem 3.2]{AS}). Hence the assumption 
\eqref{nnnnnnn.1} of Theorem~\ref{nnnnnnn} is satisfied. 

From \cite[(2.3) and Theorem~2.3]{AS} we have 
\begin{equation}\label{eqeq101}
F(P(\lambda))=F(\Delta(\lambda))\cong\nabla(\lambda)\cong I(\lambda),
\end{equation}
if $\la$ is regular. By \cite[Theorem~3.2]{AS},
$F$ commutes with projective functors. Applying projective functors to
\eqref{eqeq101} gives $F(P(\mu))=I(\mu)$ for any $\mu\in W\cdot \la$. Hence,
the assumption~\eqref{nnnnnnn.2} of Theorem~\ref{nnnnnnn} is satisfied. It is
left to verify the assumption~\eqref{nnnnnnn.3} of Theorem~\ref{nnnnnnn}.

Since the endomorphism ring of $P(w_0\cdot\la)$ is symmetric, by Proposition~\ref{cperfect} it is left to
check that $T_{w_0}$ is isomorphic to the identity functor when restricted to the
category of projective-injective modules. By \cite[Theorem~4]{KM}, there is
a natural transformation, $T_{w_0}\to \ID$, which is an isomorphism, when
restricted to projective-injective modules (\cite[Proposition~5.4]{AS}). In
particular, the assumption~\eqref{nnnnnnn.3} of Theorem~\ref{nnnnnnn} is satisfied. 
Theorem~\ref{nnnnnnn} therefore implies that $\cL(T_{w_0})^2$ is a Serre functor of
$\mathcal{D}^b(\cO_\la)$. The first part of the proposition follows. 

Let now $\la$ be dominant, integral and regular.  We again want to apply
Theorem~\ref{nnnnnnn}. The functor
$F=(C_{w_0})^2:\cO_\la\rightarrow\cO_\la$ is both right exact and 
additive by definition. Its derived functor defines a self-equivalence of
$\mathcal{D}^b(\cO_\la)$ by \cite[Theorem~5.7]{MS}. That
$F(P(\mu))=I(\mu)$ for any $\mu\in W\cdot \la$ follows inductively from
\cite[Proposition~3.1]{Irshuffle}, \cite[Theorem~5.7, Lemma~5.2 and Proposition~5.3]{MS}.
Since $\END_\mg(P(w_0\cdot\la))$ is symmetric, it is, by
Proposition~\ref{cperfect}, left to
check that $F$ is isomorphic to the identity functor when restricted to the
category of projective-injective modules. That $F$ preserves projective-injective
modules follows from \cite[Theorem~4.1(1)]{Irshuffle}. 
From \cite[Section~2.4]{Sperv} it follows that $F$ commutes with 
the action of the centre of $A$, which, because of the double centraliser
and commutativity of $\END_\mg(P(w_0\cdot\la))$, is in fact $\END_\mg(P(w_0\cdot\la))$. This
implies that $F$, restricted to the category of projective-injective modules,
is isomorphic to the identity functor. Theorem~\ref{nnnnnnn} now implies that 
$\cL(C_{w_0})^2$ is a Serre functor of $\mathcal{D}^b(\cO_\la)$. From the
uniqueness of Serre functors we get in particular $\cL\,(T_{w_0})^2\cong\cL\,(C_{w_0})^2$. 
\end{proof}

We obtain the following surprising consequence:

\begin{cor}
\label{Arkhipov}
Let $\la$ be an integral, dominant and regular weight. Considered as endofunctors of
$\cO_\la$, there is an
 isomorphism of functors $(C_{w_0})^2\cong (T_{w_0})^2$. In particular
$(C_{w_0})^2$ commutes with projective functors.   
\end{cor}

\begin{proof} 
The functors are isomorphic when restricted to the additive category of 
projective modules, since they
both give rise to a Serre functor. On the other hand, they are both right
exact and $\cO_\la$ has finite global dimension. Therefore, the isomorphism
extends uniquely to the whole category $\cO_\la$. Twisting functors commute with
projective functors (see \cite[Section~3]{AS}), hence $(C_{w_0})^2$ commutes with
projective functors as well. 
\end{proof}

\begin{remark}
\label{commute}
{\rm We would like to draw the reader's attention to the following
  observations concerning the {\it principal block $\cO_0$} of $\cO$:
\begin{enumerate}[1.)]
\item 
The functor $\op{d}T_{w_0}^2\op{d}$ is exactly Enright's completion functor, see \eg \cite{Jo1}. This follows from \cite[Section~3]{KM}.
\item Considered as an endofunctor of $\cO_0$, the functor $C_{w_0}$ does not commute with the action of the centre 
of the universal enveloping algebra of $\mathfrak{g}$ (or with the centre
of $\cO_0$) and does not commute with translation functors even if 
$\mg=\mathfrak{sl}_2$ (whereas $T_{w_0}$ does, see \cite[Section~3]{AS}). 
This is because $C_{w_0}$ twists the action of the centre by $w_0$ (this follows 
from \cite[Section~2.4]{Sperv}). This means, however, that $C_{y^{-1}} C_y$ 
commutes with the action of the centre of the category for any $y\in W$
(however, not necessarily with projective functors).
\item Since $(C_{w_0})^2$ induces the identity on the category of injective
modules, it follows that $(C_{w_0})^4\cong (C_{w_0})^2$. It
is easy to see that already $T_{w_0}^3\cong T_{w_0}^2$. Moreover,  the functors  
$C_{w_0}$, $(C_{w_0})^2$, and  $C_{w_0}^3$  are pairwise non-isomorphic; and  
the functors  $T_{w_0}$, $T_{w_0}^2$ are not isomorphic. 
\item If $w_0=s_1\dots s_k$ is a reduced decomposition, then 
\begin{gather*}
\begin{array}{rcl}
\mathcal{L}(C_{w_0})^2&\cong&\left(\mathcal{L}C_{s_1}\,\mathcal{L}C_{s_2}\cdots
\mathcal{L}C_{s_t}\right)\left(\mathcal{L}C_{s_1}\,\mathcal{L}C_{s_2}\cdots
\mathcal{L}C_{s_t}\right),\\
\mathcal{L}(T_{w_0})^2&\cong&\left(\mathcal{L}T_{s_1}\,\mathcal{L}T_{s_2}\cdots
\mathcal{L}T_{s_t}\right)\left(\mathcal{L}T_{s_1}\,\mathcal{L}T_{s_2}\cdots
\mathcal{L}T_{s_t}\right).
\end{array}
\end{gather*}
The first isomorphism follows for example from \cite[Proposition~3.1]{Irshuffle} 
by standard arguments. The second follows directly from
\cite[Theorem~2.2, Theorem~2.3 and Theorem~3.2]{AS}.
\item The Serre functor $\mathbb{S}$ for $\cD^b(\cO_0)$  satisfies
$\mathbb{S}^k\not\cong \mathbb{S}^l$ for all $k\neq l$. Indeed,
from \cite[Corollary~6.2]{AS} it follows that $\mathbb{S}^kL(0)\cong L(0)[k2 l(w_0)]$.
From \cite[Corollary~6.2]{AS} it also follows that $\mathbb{S}^k \not\cong[l]$ for any $k,l$ 
because $\mathbb{S}P(w_0\cdot 0)\cong P(w_0\cdot 0)$.
\item The braid group acts on $\cD^b(\cO_0)$ via the auto-equivalences
$\mathcal{L}C_{s}$ and via the auto-equivalences $\mathcal{L}T_{s}$. Since the Serre
functor commutes with auto-equivalences, it is natural to expect that it should 
correspond to a central element in the Braid group. In fact, $s_1\dots s_ks_1\dots s_k$ 
(see notation above) generates the centre of the Braid group $B_n$, $n\geq 3$,
see for example \cite[Corollary~1.8.4]{Birman}.
\end{enumerate}
}
\end{remark}

\subsection{The parabolic category $\cO$ in the sense of Rocha-Caridi}
\label{parabolic}

Our next task is to describe the Serre functor for the bounded derived
category associated with the principal block of a parabolic category $\cO$ in
the sense of \cite{RC}. The situation here is much more complicated, since
there are in general non-isomorphic indecomposable projective-injective modules in the same
block and we do not yet know if the endomorphism ring of a basic
projective-injective module is symmetric. However, the knowledge of the
Serre functor for the bounded derived category of $\cO_0$ turns out to be
extremely useful to determine the Serre functor for the parabolic situation. \\
  
Let $\p\supset\mathfrak{b}$ be a parabolic subalgebra of 
$\mg$ with corresponding Weyl group $W_\p\subset W$. Let $w_0^\p$ be the longest 
element in $W_\p$. For any integral dominant weight $\la$ let $\cO_\la^\p$ be the full subcategory of $\cO_\la$ 
given by locally $\p$-finite objects. This category was introduced in \cite{RC}. 
For any $w\in W$ let $\Delta^\p(w\cdot\la)$ denote the corresponding parabolic Verma
module with highest weight $w\cdot\la$, \ie the maximal quotient,
contained in $\cO_\la^\p$, of the Verma module $\Delta(w\cdot\la)\in\cO_\la$. 
Note that $\Delta^\p(w\cdot\la)\not=0$ if and only if $w$ is a shortest coset
representative in $W_\p\backslash W$.   

Let from now on $\la$ be dominant, integral and regular. Since any object in $\cD^b(\cO_\la)$ for which all cohomology objects
are contained in $\cO_\la^\p$ is quasi-isomorphic to some complex of objects
from $\cO_\la^\p$ (see the proof of \cite[Proposition~1.7.11]{KaSha}), 
$\cD^b(\cO_\la^\p)$ embeds as a full triangulated
subcategory in $\cD^b(\cO_\la)$. Note that translations through walls preserve
the parabolic subcategory. We may therefore consider the restriction
of $\cL C_{w_0}$ to this subcategory as well as to the subcategory
$\cO_\la^\p$ (considered as a subcategory of $\cD^b(\cO_\la^\p)$). 
We get the following result: 

\begin{prop}
\label{Serreparab}
For any integral dominant and regular weight $\la$ we have: 
\begin{enumerate}[(1)]
\item  \label{first} The functor $\cL(C_{w_0})[-l(w_0^\p)]$ maps parabolic Verma modules to
  parabolic dual Verma modules. More precisely 
\begin{equation}
\label{Vermas}  
\cL(C_{w_0})[-l(w_0^\p)]\,\Delta^\p(w\cdot\la)\cong\op{d}\Delta^\p(w_0^\p ww_0\cdot\la)
\end{equation}
for any parabolic Verma module $\Delta^\p(w\cdot\la)\in\cO_\la^\p$.
\item The category $\cO_\la^\p$ is Ringel self-dual. 
\item \label{secondlast} The functor $\cL(C_{w_0})^2[-2l(w_0^\p)]$ 
maps projectives in $\cO^\p$ to injectives in $\cO^\p$.
\item \label{last} The functor $\cL(C_{w_0})^2[-2l(w_0^\p)]$ is a Serre
functor for $\cD^b(\cO_\la^\p)$. 
\end{enumerate}
\end{prop}

\begin{proof} 
We first check that $w_0^\p ww_0$ is indeed a shortest coset representative, if
so is $w$. Let $s\in W_\p$ be a simple reflection. Then
\begin{multline*}
l(sw_0^\p w w_0)=l(w_0)-l(sw_0^\p w)=l(w_0)-l(sw_0^\p)-l(w)=\\
=l(w_0)-l(w_0^\p)+1-l(w)= l(w_0)-l(w_0^\p w)+1=l(w_0^\p w w_0)+1.
\end{multline*}

Let $\Delta^\p(w\cdot\la)$ be a parabolic Verma module in $\cO_\la^\p$. From
\cite[Section 4]{L2} we have a finite resolution, $P^w_\bullet$, of 
$\Delta^\p(w\cdot\la)$ by Verma modules, where
\begin{eqnarray*}
P^w_i=\displaystyle\bigoplus_{y\in W_\p, l(y)=i} \Delta(yw\cdot\la).
\end{eqnarray*}
The involved maps are clear up to scalars, for the exact normalisation we
refer to \cite[Lemma 4.1]{L2}. 

For a simple reflection, $s$, the definition of $C_s$ implies 
$C_s\Delta(x\cdot \la)\cong \Delta(xs\cdot \la)$ if $xs>x$, and 
$C_s\nabla(x\cdot \la)\cong \nabla(xs\cdot \la)$ if $xs<x$. Therefore
$C_{w_0}\Delta(x\cdot \la)\cong \nabla(xw_0\cdot \la)$, which implies that
$C_{w_0}$ is exact on the category of modules with Verma flag. This gives
$\cR C_{w_0}\Delta^\p(w\cdot\la)\cong C_{w_0}(P^w_\bullet)$. 
Dually, $\op{d}P^{w_0^\p ww_0}_\bullet$ is a coresolution of 
$\op{d}\Delta^\p(w_0^\p ww_0\cdot\la)$. On the other hand, applying 
$C_{w_0}$ to the resolution $P^{w}_\bullet$ gives a complex, $Q_\bullet$, where 
\begin{displaymath}
Q_i=\displaystyle\bigoplus_{y\in W_\p, l(y)=i}
\op{d}\Delta(yww_0\cdot\la)=
\displaystyle\bigoplus_{y\in W_\p, l(y)=i}
\nabla(yww_0\cdot\la) 
\end{displaymath}           
The maps in this complex satisfy the dual version of
\cite[Lemma 4.1]{L2}. Hence $\op{d}P^{w_0^\p ww_0}_\bullet\cong Q_\bullet[-l(w_0^\p)]$. The 
formula \eqref{Vermas} follows. 

Let now $F=\cL_{l(w_0^p)}C_{w_0}$. We claim that $F$, restricted to $\cO_\la^\p$, is right exact. Note that the formulas above imply that
$\cL_i C_{w_0}M=0$ for any $M\in\cF(\Delta^{\p})$ and $i\not=l(w_0^\p)$,
in particular, for $M\in\cO_\la^\p$ projective. Let $M\in\cO_\la^\p$ be
arbitrary. Choose a short exact sequence $K\hookrightarrow
P\twoheadrightarrow M$, where
$P\in\cO_\la^\p$ is projective. Since the global dimension of $A_\la^\p$ is finite,
one obtains $\cL_i C_{w_0}M=0$ for all $M$ and all $i<l(w_0^\p)$ by induction.
Therefore, $F$ is right exact. It is known (see \eg 
\cite[Lemma~5.1,~Lemma~5.2]{MS}) that
$G=\cR^{l(w_0^\p)}(\op{d}C_{w_0}\op{d})$ is the right adjoint functor of
$F$. From the formula~\eqref{Vermas} it follows that $F$ defines an
equivalence $\cF(\Delta^{\p})\cong\cF(\nabla^{\p})$ with
inverse $G$. Proposition~\ref{Ringel2} implies that $\cO_\la^\p$ is Ringel
self-dual. Applying Proposition~\ref{Ringel2} twice, we get that the functor
$\cL(C_{w_0})^2[-2l(w_0^\p)]$ maps projective modules to injective modules. To prove
that $\cL(C_{w_0})^2[-2l(w_0^\p)]$ is a Serre functor we only have to verify
the last assumption of Theorem~\ref{nnnnnnn}. This is not completely
trivial. Instead of applying again Theorem~\ref{nnnnnnn} we will give an
alternative argument after the following lemma.
\end{proof}

Let $\la$ still be dominant, integral and regular. Let $\mathfrak{i}:\cO_\la^\p\to \cO_\la$ denote the exact inclusion functor,
let $\mathrm{Z}:\cO_\la\to \cO_\la^\p$ be its left adjoint and
$\hat{\mathrm{Z}}:\cO_\la\to \cO_\la^\p$ be the right adjoint to
$\mathfrak{i}$ (\ie $\mathrm{Z}$ is the {\em Zuckerman} functor of taking
the maximal quotient in $\cO_\la^\p$, $\hat{\mathrm{Z}}\cong
\mathrm{d}\,\mathrm{Z}\,\mathrm{d}$). To proceed we will need the following
result from folklore:

\begin{lemma}\label{Zucker}
There are isomorphisms of functors: 
\begin{eqnarray*}
\mathrm{d}\,\mathfrak{i}\,\mathcal{L}\mathrm{Z}\,\mathrm{d}\cong
\mathfrak{i}\,\mathcal{L}\mathrm{Z}[-2l(w_0^\p)]:&&\cD^b(\cO_\la)\rightarrow\cD^b(\cO_\la)\\  
\mathrm{d}\,\,\mathcal{L}\mathrm{Z}\,\mathrm{d}\cong\mathcal{L}\mathrm{Z}[-2l(w_0^\p)]:&&\cD^b(\cO_\la)\rightarrow\cD^b(\cO_\la^\p) 
\end{eqnarray*}
\end{lemma}

\begin{proof}
Using \cite[Proposition~4.2]{EW} we can fix an isomorphism,
\begin{displaymath}
\mathrm{d}\,\mathfrak{i}\,\mathcal{L}\mathrm{Z}\,\mathrm{d}\Delta(\lambda)\cong
\mathfrak{i}\,\mathcal{L}\mathrm{Z}\Delta(\lambda)[-2l(w_0^\p)].
\end{displaymath}
Since  $\mathrm{Z}$ commutes with projective functors (see \eg \cite[Proposition3]{BFK}), this isomorphism lifts
to an isomorphism on projective modules. We have to
verify that it is functorial. Without loss of generality we may
assume that $\p=\p_s$ is the parabolic subalgebra corresponding to a simple
reflection $s$. The general case for arbitrary $\p$ follows then by
induction. Associated to $s$, there
is a complex of functors 
\begin{eqnarray*}
  \mathrm{T}_s\rightarrow \ID\rightarrow\mathfrak{i}\,\mathrm{Z} 
\end{eqnarray*}
which gives rise to a short exact sequence when applied to projective objects
(\cite[Theorem 4]{KM} and \cite[Proposition~5.4]{AS}). Taking the left
derived functors we therefore get an isomorphism
$\cL_2(\mathfrak{i} \mathrm{Z})\cong\cL_1\mathrm{T}_s$ (see \cite[Proposition
1.8.8]{KaSha}). From \cite[Theorem~4.1]{AS} and \cite[Theorem 1]{MazStr} we
have an isomorphism of functors
$\cL_1\mathrm{T}_s\cong\mathrm{d}\,\mathfrak{i}\,\mathrm{Z}\,\mathrm{d}$. The
first statement follows, the second is then also clear.  
\end{proof}

Now we are ready to complete the proof of Proposition~\ref{Serreparab}:
\begin{proof}[Proof of the last part of Proposition~\ref{Serreparab}]
To prove that $\cL(C_{w_0})^2[-2l(w_0^\p)]$ is a Serre
functor for $\cD^b(\cO_\la^\p)$ it is enough to show that for the
functor $\mathrm{G}=\mathrm{Z}\,\mathcal{L}_{2l(w_0^\p)}(C_{w_0})^2\,\mathfrak{i}$
we have an isomorphisms, natural in
both arguments, as follows:
\begin{equation}\label{eq203}
\mathrm{Hom}_{\cO_\la^\p}(P^\p,\mathrm{G} P^\p)\cong
\mathrm{Hom}_{\cO_\la^\p}(P^\p,P^\p)^*,
\end{equation}
where $P^\p$ is a projective generator of $\cO_\la^\p$. Without loss of generality we assume $P^\p=\mathrm{Z}P$, 
where $P$ is a projective generator of $\cO_\la$. We have isomorphisms
\begin{displaymath}
\begin{array}{rcl}
\mathrm{Hom}_{\cO_\la^\p}(P^\p,\mathrm{G} P^\p) & \cong & 
\mathrm{Hom}_{\cO_\la^\p}(\mathrm{Z}P,\mathrm{G} P^\p)  \\
& \cong & \mathrm{Hom}_{\cO_\la}(\mathfrak{i}\,\mathrm{Z}P,\mathfrak{i}\,\mathrm{G} P^\p)  \\
& \cong & \mathrm{Hom}_{\cO_\la}(P, \mathcal{L}_{2l(w_0^\p)}(C_{w_0})^2\,\mathfrak{i}P^\p)  \\
\text{(by Proposition~\ref{Serreparab}~\eqref{first})}
& \cong & \mathrm{Hom}_{\cD^b(\cO_\la)}(P, \mathcal{L}(C_{w_0})^2\,\mathfrak{i}P^\p[-2l(w_0^\p)]) \\
& \cong & \mathrm{Hom}_{\cD^b(\cO_\la)}(P[2l(w_0^\p)], \mathcal{L}(C_{w_0})^2\mathfrak{i}P^\p)  \\
\text{(by Proposition~\ref{prcato})}
& \cong & \mathrm{Hom}_{\cD^b(\cO_\la)}(\mathfrak{i}P^\p,P[2l(w_0^\p)])^* \\ 
\text{(by adjointness of $\mathfrak{i}$ and $\mathrm{d}\mathrm{Z}\mathrm{d}$)}
& \cong & \mathrm{Hom}_{\cD^b(\cO_\la^\p)}(P^\p,\mathrm{d}\,\mathcal{L}\mathrm{Z}\,\mathrm{d}P[2l(w_0^\p)])^* \\
\text{(by Lemma~\ref{Zucker})}
& \cong & \mathrm{Hom}_{\cD^b(\cO_\la^\p)}(P^\p,\mathcal{L}\mathrm{Z}P)^* \\
& \cong & \mathrm{Hom}_{\cO_\la^\p}(P^\p,\mathrm{Z}P)^*  \\
& \cong & \mathrm{Hom}_{\cO_\la^\p}(P^\p,P^\p)^*,
\end{array}
\end{displaymath}
which are natural in both arguments. This completes the proof of Proposition~\ref{Serreparab}.
\end{proof}

As an application we get the following nontrivial result
\begin{theorem}
\label{symmetry}
  Let $\la$ be an integral, regular and dominant weight. We consider the
  category $\cO_\la^\p$, where  
  $\p\supset\mathfrak{b}$ is some parabolic subalgebra of $\mg$. Let $Q$ be a
  basic projective-injective module in $\cO_\la^\p$. Then
  $\END_\mg(Q)$ is symmetric.  
\end{theorem}

For the proof we need the following 

\begin{lemma}
\label{esspartdom}
In the situation of Theorem~\ref{symmetry} we have the following: the socle $S$ of
  $\Delta^\p(\la)$ is simple and $\HOM_\mg(Q,\Delta^\p(\la)/S)=0$. 
\end{lemma}

\begin{proof}
  Assume $L$ is a composition factor of $\Delta^\p(\la)$ such that
  $\HOM_\mg(Q,L)\not=0$. From the latter it follows that the projective cover,
  $P$, of $L$ is projective-injective, hence tilting. Since $L$ is a composition factor of
  $\Delta^\p(\la)$, we have $\HOM_\mg(\Delta^\p(\la),P)\not=0$. 
Therefore, $L(\la)$ appears as a composition factor in $P$. Hence $P=T^\p(\la)$
  and $L$ is unique. On the other hand
\begin{displaymath}
1=[T^\p(\la)\::\:\Delta^\p(\la)]=[P\::\:\Delta^\p(\la)]=[\Delta^\p(\la)\::\:L].
\end{displaymath}
  Since (by \cite{Irself}) any simple module appearing in the socle of a parabolic Verma module
  is not annihilated by $\HOM_\mg(Q,\bullet)$, the statement follows.  
\end{proof}

\begin{proof}[Proof of Theorem~\ref{symmetry}]
  By Example~\ref{examples}~\eqref{exparabO} and Proposition~\ref{cperfect} 
  it is enough to prove that the Serre functor is
  isomorphic to the identity functor when restricted to the additive subcategory
  of projective-injective modules. Let $\mS=\cL(C_{w_0})^2[-2l(w_0^\p)]$ be
  the Serre functor of $\cD^b(\cO_\la^\p)$. The idea of the proof is the following:  From
  Proposition~\ref{Serreparab} and Corollary~\ref{Arkhipov} we know that
  $S:=\mathfrak{i}\,\cL_{2l(w_0^\p)}(C_{w_0})^2\, \mathrm{Z}:\cO_\la\rightarrow\cO_\la$ is right exact and commutes with translations through walls, even in a natural
  way as defined in \cite{projfunc}. We will construct
  another right exact functor $\mathrm{G}:\cO_\la\rightarrow\cO_\la$ which again
  naturally commutes with translations through walls and coincides with $S$
  when evaluated at $\Delta(\la)$. The main result of \cite{projfunc} states
  that two right exact additive functors, $F_1$,
  $F_2:\cO_\la\rightarrow\cO_\la$,  which agree on
  $\Delta(\la)$, and both
  naturally commute with translations through walls, are in fact
  isomorphic. From this fact we will deduce an
  isomorphism of functors $S\cong \mathrm{G}$. Even more, the main
  step is to show that we can choose a functor $\mathrm{G}$ with the additional 
  property that $\mathrm{G}=\mathfrak{i}\,\mathrm{G}'\,\mathrm{Z}$ for some $\mathrm{G}':\cO_\la^\p\rightarrow\cO_\la^\p$, and $\mathrm{G}'$ is isomorphic to the
  identity functor when restricted to the category of projective-injective
  modules in $\cO_\la^\p$. Since
  $\mathrm{Z}$ is dense and full and  $S=\mathfrak{i}\,\cL_{2l(w_0^\p)}(C_{w_0})^2\, \mathrm{Z}\cong \mathfrak{i}\,\mathrm{G}'\,\mathrm{Z}$, we get that $\mS$
  must be isomorphic to the identity functor on the subcategory of $\cO_\la^\p$ formed by projective-injective
  modules. This will finally imply the assertion of the theorem. 

Let's do the work! We have $\mS=\cL(C_{w_0})^2[-2l(w_0^\p)]$, the Serre
functor of $\cD^b(\cO_\la^\p)$, given by restriction of
$\cL(C_{w_0})^2[-2l(w_0^\p)]:\cD^b(\cO_\la)\rightarrow\cD^b(\cO_\la)$. Put
$S:=\mathfrak{i}\,\cL_{2l(w_0^\p)}(C_{w_0})^2\, \mathrm{Z}$, considered as a functor
$\cO_\la\rightarrow\cO_\la$. This functor is clearly right exact and
additive. Since $\la$
is regular, the category $\cO_\la$ is a category with full projective functors
(\cite[Proposition 16]{projfunc})
in the sense of \cite[Section 2]{projfunc}, where the projective functors are
given by compositions of translations through walls and their direct
summands. Recall from \cite[Definition 2]{projfunc} that a functor
$G:\cO_\la\rightarrow\cO_\la$ {\it naturally commutes with projective
  functors} if for any projective functor $\theta$, there is an isomorphism
of functors $\varphi_\theta:\theta\,G\cong G\,\theta$ such that the following holds:
for any two projective functors $\theta_1$, $\theta_2$ and any natural
transformation $\alpha\in\HOM(\theta_1,\theta_2)$ the following diagram commutes:
 \begin{eqnarray*}
\xymatrix{
\theta_1 G\ar@{->}[rr]^{\alpha_G}\ar@{->}[d]^{\varphi_{\theta_1}}&&\theta_2 G\ar@{->}[d]_{\varphi_{\theta_2}}\\
G\,\theta_1\ar[rr]^{G(\alpha)}&&G\,\theta_2.
}
\end{eqnarray*}
(Note the typos in the original formulation \cite[Definition 2]{projfunc}.)\\

{\it Claim 1: The functor $S:\cO_\la\rightarrow\cO_\la$ naturally commutes with projective functors.}

\begin{proof}[Proof of Claim 1]
To see this note first that $S\cong \cL_{2l(w_0^\p)}
(C_{w_0})^2\,\mathfrak{i}\,\mathrm{Z}$, where $\cL
(C_{w_0})^2:\cD^b(\cO_\la)\rightarrow\cD^b(\cO_\la)$. (This is clear
from the definitions.) From \cite[Section 6.2]{projfunc} we know that
$\mathfrak{i}\,\mathrm{Z}$ naturally commutes with projective functors on
$\cO_\la$. In \cite[Section 6.5]{projfunc} it is proved that twisting functors
on $\cO_\la$ naturally commute with projective functors on $\cO_\la$. Corollary~\ref{Arkhipov} gives an isomorphism $(T_{w_0})^2\cong
(C_{w_0})^2:\cO_\la\rightarrow\cO_\la$, hence  $(C_{w_0})^2$ naturally commutes with projective functors. From \cite[Lemma 8]{projfunc} it follows that $\cL_{2l(w_0^\p)}(C_{w_0})^2$ and $\cL_{2l(w_0^\p)}(C_{w_0})^2$ naturally commute with projective functors, therefore
so does $S$, because it is a composition of functors which commute naturally 
with projective functors. This implies Claim 1. 
\end{proof}

Let now $J:\cO_\la\rightarrow\cO_\la$ be the partial coapproximation with
respect to $M$, where $M=\oplus_{x\in W'} P(x\cdot\la)$ and
\begin{equation*}
  W'=\{x\in W\mid \mathrm{Z} P(x\cdot\la)=0 \text{ or  }
  \mathrm{Z}P(x\cdot\la)\in\cO_\la^\p \text{ is projective-injective}\}.
\end{equation*}
Recall that, when restricted to projective objects, $J$ is nothing else than
taking the trace with respect of $M$. The functor $J$ is additive and right
exact.\\

{\it Claim 2: The functor $J:\cO_\la\rightarrow\cO_\la$ naturally commutes with projective functors.}

\begin{proof}[Proof of Claim 2]
Let $\theta:\cO_\la\rightarrow\cO_\la$ be a projective functor. We first show
that $\theta\op{Tr}_M P=\op{Tr}_M\theta P$, via the natural inclusions 
\begin{eqnarray*}
  \theta\op{Tr}_M P\hookrightarrow\theta P\hookleftarrow\op{Tr}_M \theta P, 
\end{eqnarray*}
for any projective module $P\in\cO_\la$. To see this consider the short exact
sequence $\op{Tr}_M P\hookrightarrow P\twoheadrightarrow N$,
where $N$ is the canonical quotient, in particular $\HOM_\mg(M,N)=0$. We claim
that $\HOM_\mg(M,\theta N)=0$. Let $\theta'$ be the adjoint functor of
$\theta$. This is of course again a projective functor and therefore we have
the following: If $\mathrm{Z}P(x\cdot\la)=0$ then
$0=\theta'\mathrm{Z}P(x\cdot\la)\cong\mathrm{Z}\theta'P(x\cdot\la)$. If
$\mathrm{Z}P(x\cdot\la)\not=0$, but $x\in W'$, then $\mathrm{Z}P(x\cdot\la)$
is projective-injective in $\cO_\la^\p$, hence so is
$\theta'\mathrm{Z}P(x\cdot\la)\cong\mathrm{Z}\theta' P(x\cdot\la)$. In particular,  
\begin{eqnarray*}
  \HOM_\mg(M,\theta N)&\cong&\HOM_\mg(\theta'M, N)\\
&\hookrightarrow&\HOM_\mg(M^n, N)\quad \text{  (for some positive integer $n$)}\\
&=&0.
\end{eqnarray*}
The definition of the trace implies that the projective cover of $\op{Tr}_MP$ is a direct
summand of some $M^n$, $n\in\mZ_{>0}$. From the arguments above it follows that
the projective cover of $\theta\op{Tr}_M P$ is also a direct summand of some
$M^n$, $n\in\mZ_{>0}$. Altogether, $\theta\op{Tr}_M P=\op{Tr}_M\theta P$
via the natural inclusions. In other words, we may fix an isomorphism
of functors $\varphi_\theta:\theta J\cong J\theta$, restricted to the
category of projective modules, such that
\begin{eqnarray}
  \label{eq:diag1}
j_\theta\circ \varphi_\theta=\theta(j):\theta\,J\rightarrow\theta,  
\end{eqnarray}
where $j:J\rightarrow\ID$ is the obvious natural transformation. In
particular, $J$ commutes with projective functors.   
We claim that this is already enough to show that $J$ {\it naturally} commutes
with projective functors. We have to check this directly using the original
definition \cite[Definition~2]{projfunc}: Let $\theta_1$,
$\theta_2:\cO_\la\rightarrow\cO_\la$ be two projective functors and let
$\alpha\in\HOM(\theta_1,\theta_2)$ be a natural transformation between
them. Consider the following diagram of functors restricted to the additive
category formed by all projective objects:
\begin{eqnarray*}
\xymatrix{\theta_1\ar[rr]^{\alpha}\ar@/_2pc/@{=}[ddd]_{\op{id}}\ar@{<-^{)}}[d]^{\theta_1(j)}&&\theta_2\ar@/^2pc/@{=}[ddd]^{\op{id}}\ar@{<-_{)}}[d]_{\theta_2(j)}\\
\theta_1J\ar[rr]^{\alpha_J}\ar@{->}[d]^{\varphi_{\theta_1}}&&\theta_2J\ar@{->}[d]_{\varphi_{\theta_2}}\\
J\theta_1\ar[rr]^{J(\alpha)}&&J\theta_2\\
\theta_1\ar[rr]^{\alpha}\ar@{<-_{)}}[u]_{j_{\theta_1}}&&\theta_2\ar@{<-^{)}}[u]^{j_{\theta_2}}.
}  
\end{eqnarray*}
The two ``squares'', the one on the left hand side and the on the right side, commute because of
\eqref{eq:diag1}. The squares at the top and bottom commute by definition (of a natural
transformation). We only have to show that the middle square commutes as well,
\ie $J(\alpha)\circ\varphi_{\theta_1}=\varphi_{\theta_2}\circ\alpha_J$. Since
$j_{\theta_2}$ is injective (on projective modules) it is enough to show that 
$j_{\theta_2}\circ
J(\alpha)\circ\varphi_{\theta_1}=j_{\theta_2}\circ\varphi_{\theta_2}\circ\alpha_J$.
Since all the other parts of the diagram commute we can calculate
\begin{eqnarray*}
  j_{\theta_2}\circ
J(\alpha)\circ\varphi_{\theta_1}&=&\alpha\circ j_{\theta_1}\circ\varphi_{\theta_1}\\
&=&\alpha\circ\theta_1(j)\\
&=&\theta_2(j)\circ\alpha_J\\
&=&j_{\theta_2}\circ\varphi_{\theta_2}\circ\alpha_J.
\end{eqnarray*}
Hence, $J$ commutes naturally with projective functors when restricted to
projective objects. Since the involved functors are right exact, Claim 2 follows.
\end{proof}
{\it Claim 3: There is an isomorphism of modules
 $S\Delta(\la)\cong\mathfrak{i}\,\mathrm{Z}\,J\,\mathfrak{i}\,\mathrm{Z}\,J\Delta(\la)$}.
\begin{proof}[Proof of Claim 3]
We first show that
$\mathrm{Z}\,J\Delta(\la)\cong\op{soc}\Delta^\p(\la)$, the socle of $\Delta^\p(\la)$. Define $U$ to be the
module which fits into the canonical short exact sequence 
\begin{eqnarray}
\label{ses}
  0\rightarrow
  U\longrightarrow\Delta(\la)\longrightarrow\Delta^p(\la)\rightarrow 0.
\end{eqnarray}
From Lemma~\ref{esspartdom} we have
$\mathrm{Z}\,J\Delta^\p(\la)\cong\op{soc}\Delta^\p(\la)$. On the other hand,
$ZU=0$ by definition, \ie the projective cover $P_U$ of $U$ is annihilated by
$\mathrm{Z}$, hence also $\mathrm{Z}\,J P_U=0$. Since
$\mathfrak{i}\,\mathrm{Z}\,J$ is right exact, the sequence~\eqref{ses} implies
$\mathrm{Z}\,J U=0$ and $\mathrm{Z}\,J\Delta(\la)\cong\op{soc}\Delta^\p(\la)$. From the double
centraliser property we have an exact sequence in $\cO_\la^\p$ of the form 
\begin{eqnarray*}
  0\rightarrow\Delta(\la)^\p\longrightarrow Q_1\longrightarrow Q_2,
\end{eqnarray*}
where $Q_1$ and $Q_2$ are projective-injective modules in $\cO_\la^\p$. 
There is therefore also an exact sequence in $\cO_\la$ of the form 
\begin{eqnarray*}
Q_2\rightarrow Q_1\rightarrow \op{d}\Delta^\p(\la)\rightarrow 0.  
\end{eqnarray*}
Note that $J\,Q_1=Q_1$, $J\,Q_2=Q_2$. Hence
$\mathfrak{i}\,\mathrm{Z}\,J\op{d}\Delta^\p(\la)\cong\op{d}\Delta^\p(\la)$. Altogether we
have
$\mathfrak{i}\,\mathrm{Z}\,J\,\mathfrak{i}\,\mathrm{Z}\,J\Delta^\p(\la)\cong\op{d}\Delta^\p(\la)$.
The latter is isomorphic to $S\Delta^\p(\la)$ by
Proposition~\ref{Serreparab}. This proves Claim 3.
\end{proof}

From \cite[Theorem 1]{projfunc} we get the
existence of an isomorphism 
\begin{eqnarray*}
  \alpha:\mathfrak{i}\,\mathrm{Z}\,J\,\mathfrak{i}\,\mathrm{Z}\,J\cong S.
\end{eqnarray*}
By definition, $J$ is isomorphic to the identity functor when restricted to
projective-injectives in $\cO_\la^\p$. Therefore, $\alpha$ induces an isomorphism
of functors  
\begin{eqnarray*}
\ID&\cong& \cL (C_{w_0})^2[-2l(w_0^\p)]
\end{eqnarray*}
when restricted to the category of projective-injective modules in
$\cO_\la^\p$. This is exactly the statement that the Serre
functor for $\cD^b(\cO_\la^\p)$ is isomorphic to the identity when restricted
to the additive category formed by projective-injective modules. The assertion
of the theorem follows then finally from Proposition~\ref{cperfect}.  
\end{proof}

We get the following consequence
\begin{cor}
  In the situation of Theorem~\ref{symmetry} the Serre functor
  for $\cD^b(\cO_\la^\p)$ is isomorphic to $\cL((\CA_Q)^2)$
\end{cor}

\begin{proof}
  This follows directly from Theorem~\ref{symmetry} and Theorem~\ref{Serrecoapprox}.
\end{proof}

\subsection{Harish-Chandra bimodules}

Let $^{}_\la\cH_\mu$ be the category of Harish-Chandra bimodules as in Examples~\ref{examples}\eqref{exHC}, where $\la$ 
and $\mu$ are integral and dominant. Recall the subcategory $^{}_\la\cH_\mu\cong A_\la^\mu\MOF$.
In \cite{BG} (see also \cite[6.17, 6.23]{Ja2}) it is proved that
$^{}_\la\cH_\mu^1$ is equivalent to the full subcategory $\mathcal{C}_\la^\mu$  of 
$\cO_\la$ given by all modules $M$, which have an exact presentation, 
\begin{equation}\label{eqpresent}
P_1\rightarrow P_2\rightarrow M\rightarrow 0,
\end{equation}
where $P_1$ and $P_2$ are projective and the simple modules
in their heads are of the form $L(x\cdot\la)$, where $x$ is a longest coset
representative in $W_\mu\backslash W /W_\la$. Note that this category does not
have finite global dimension in general. Nevertheless we have the following 

\begin{prop}
\label{SerreHC}
Let $\la$, $\mu$ be integral dominant weights. Then 
\begin{enumerate}[(1)]
\item $A_\la^\mu$ is Ringel self-dual. 
\item $\mathcal{D}_{\mathrm{perf}}(A_\la^\mu)$ has a Serre functor,
namely $\mathcal{L}(C_{w_0})^2$ (via the identification of
$\mathcal{C}_\la^\mu$ with $A_\la^\mu\MOF$ given by \cite[Section~5]{Au}).
\end{enumerate}
\end{prop}

\begin{proof}
From \cite[Proposition~4.2]{MS} and its dual version we get that $C_{w_0}$ maps 
standard modules to costandard modules. From \cite[Lemma 5.18]{MS} it follows 
that $C_{w_0}$ defines an equivalence between the categories of modules with standard flag and 
modules with costandard flag. In particular, $\Ext^1_{A_\la^\mu}(C_{w_0}P,\nabla)\cong
\Ext^1_{A_\la^\mu}(P,C_{w_0}^{-1}\nabla)=0$ for any costandard module $\nabla$ and any 
projective $P$. Hence $C_{w_0}$ maps a minimal projective generator to a characteristic 
cotilting module, which is tilting (see \eg \cite[Section~6]{FKM1}). 
Since it is an equivalence, it preserves the endomorphism ring. The first part of the 
theorem follows.

To prove the second part of the theorem we  use Proposition~\ref{prcato}.
Because of this result and \cite[Proposition~20.5.5(i)]{Ginzburg} it is
enough to show that $\mathcal{L}(C_{w_0})^2$ preserves 
$\mathcal{D}_{\mathrm{perf}}(A_\la^\mu)$. Since $\mathcal{L}(C_{w_0})^2$ is
a Serre functor for $\cO$, it sends indecomposable projective modules
to the corresponding indecomposable injective modules. However,
an injective $A_\la^\mu$-module has a presentation of the form
\eqref{eqpresent} by \cite[Corollary~2.11]{MS}. Since all tilting
$A_\la^\mu$-modules are also cotilting (\cite[Section~6]{FKM1}),
it follows that injective $A_\la^\mu$-modules have finite projective 
dimension. Hence  $\mathcal{L}(C_{w_0})^2$ 
preserves $\mathcal{D}_{\mathrm{perf}}(A_\la^\mu)$ and the statement follows.
\end{proof}

\begin{remark}\label{remhch}
{\rm
From \cite[Theorem~2.5]{Dlab} it follows that $A_\la^\mu$ has finite global dimension
if and only if the standard and proper standard $A_\la^\mu$-modules coincide
(i.e. $A_\la^\mu$ is  quasi-hereditary). Using the description of standard modules as 
in \cite[Proposition~2.18]{MS} it is easy to see that this is the case if and only if
$\mu$ is regular or $A_\la^\mu$ is semi-simple.
}
\end{remark}

\begin{remark}\label{remhch2}
{\rm
One can also show that for any  $N,P\in A_\la^\mu\MOF$, 
where $P$ is projective, the Serre functor $\mathcal{L}(C_{w_0})^2$ from Proposition~\ref{prcato} induces a natural isomorphism, 
\begin{displaymath}
\HOM_{A_\la^\mu}(N, (C_{w_0})^2P)\cong\HOM_{A_\la^\mu}(P,N)^\ast.
\end{displaymath}
}
\end{remark}

\subsection{The category $\cO$ for rational Cherednik algebras}

We briefly recall the facts about rational Cherednik algebras which are
important in our setup. We refer for example to \cite{GGOR} for details.   
Let $V$ be a finite dimensional vector space, $\mathtt{W}\subset \op{GL}(V)$ a
finite reflection group, and
$\mC[\mathtt{W}]$ the group algebra of $\mathtt{W}$ over $\mC$. Let $\cA$ denote the set of reflection
hyperplanes. If $h\in\cA$ then $\mathtt{W}_h$ denotes the (pointwise) stabiliser of
$h$ in $\mathtt{W}$. Let $\gamma:\cA\rightarrow \mC[\mathtt{W}]$ be a $\mathtt{W}$-equivariant map such that
$\gamma(h)\in \mC[\mathtt{W}_h]\subset \mC[\mathtt{W}]$. Associated to the pair
$(V,\gamma)$ we have $\op{H}=\op{H}(V,\gamma)$, the corresponding
rational Cherednik algebra as defined and studied for example in
\cite{GGOR}. As a vector space, $\op{H}(V,\gamma)$ is isomorphic to $S(V)\otimes \mC[\mathtt{W}]\otimes
S(V^\ast)$, where $S(V)$ denotes the algebra of polynomial functions in
$V^\ast$ (This is the PBW-theorem \cite[Theorem 1.3]{EG}). The occurring
three algebras $S(V)$, $S(V^\ast)$ and $\mC[\mathtt{W}]$ are in fact
subalgebras, for the nontrivial commutator relations between them (involving
the parameter $\gamma$) we refer to \cite{GGOR}, \cite{Guay}. 
Let $\cO=\cO(\op{H}, V,\gamma)$ be
the corresponding category $\cO$ given by all finitely generated
$\op{H}$-modules which are locally $S(V^\ast)$-finite. This is a highest weight category, where
the isomorphism classes of simple modules  are in natural bijection with irreducible modules for
$\mathtt{W}$. More precisely, if $E$ is an irreducible $\mathtt{W}$-module, then
$\Delta(E)=\op{H}(V,\gamma)\otimes_{B}E$, where $B$ is the subalgebra of
$\op{H}(V,\gamma)$ generated by
$S(V^\ast)$ and $\mC[\mathtt{W}]$. (The action of $p\in S(V^\ast)$ on $E$ is
given by multiplication with $p(0)$). The simple head $L(E)$ of $\Delta(E)$ is the simple
module corresponding to $E$. 

In general, $\op{H}(V,\gamma)$ is not isomorphic to its opposite algebra
$\op{H}(V,\gamma)^{\opp}$, therefore there is no simple preserving duality. However, we have an isomorphism (\cite[Section 4.2]{GGOR}) of
algebras $\op{H}(V,\gamma)\cong \op{H}(V^\ast,\dag\circ\gamma)^{opp}$, where $\dag:
\mathtt{W}\rightarrow \mathtt{W}, g\mapsto g^{-1}$ (the isomorphism is given by extending
$\dag$ trivially to $S(V^\star)$ and sending $v\in V$ to $-v$).

With this fixed isomorphism one can define two {\it contravariant} functors,
namely 
\begin{itemize}
\item the {\it naive duality} (\cite[Proposition 4.7]{GGOR}):
\begin{eqnarray*}
\op{d}_{V,\gamma}:\cO(V,\gamma)
&\rightarrow& \cO(V^\ast,\dag\circ\gamma)
\end{eqnarray*}
 by sending an object $M$ to the largest submodule $\op{d}_{V,\gamma}(M)$ 
of (the ordinary vector space dual) $M^\ast$ which is locally
$S(V^\ast)$-finite. This is a right $\op{H}(V,\gamma)$-module and becomes a
left 
$\op{H}(V^\ast,\dag\circ\gamma)$-module via the fixed isomorphism.
This duality sends the simple module $L(E)$ to the simple module $L(\check{E})$ indexed by
the dual representation $\check{E}$ of $E$. Projective objects are sent to injectives and
standard objects to costandard objects. 
\item the functor $D_{V,\gamma}$ (see \cite[Proposition 4.10]{GGOR}) 
  \begin{eqnarray*}
 D_{V,\gamma}=\EXT^{\DIM V}_{\op{H}(V,\gamma)}(\bullet,\op{H}(V,\gamma)):\cO(V,\gamma)
 \rightarrow\cO(V,\dag\circ\gamma). 
  \end{eqnarray*}
\end{itemize}

\begin{conj}
\label{conjCheredsym}
  Let $\op{H}(V,\gamma)$ be a rational Cherednik algebra with the corresponding
  category $\cO(V,\gamma)$. Then 
  \begin{eqnarray*}
S=\op{d}_{V^\ast,\gamma^\dag}\,
  D_{V^\ast,\gamma}\,\op{d}_{V,\dag\circ\gamma}\, D_{V,\gamma} 
  \end{eqnarray*}
is right exact and $\cL S$ is a Serre functor. 
\end{conj}

To prove this conjecture it would be enough to verify the assumptions in Theorem~\ref{nnnnnnn}, where
  $F=F_1F_2$, 
  $F_1=\op{d}_{V,\dag\circ\gamma}\, D_{V,\gamma}$ and $F_2=\op{d}_{V^\ast,\gamma^\dag}\,
  D_{V^\ast,\gamma}$. The fact that $F$ is right exact follows directly from
  \cite[4.1]{GGOR}. The assumption \eqref{nnnnnnn.1} follows directly from \cite[Lemma
  4.1, Proposition 4.7]{GGOR}. The assumption~\eqref{nnnnnnn.2} is proved in
  \cite[Proposition 5.21]{GGOR}. We do not know if
  assumption~\eqref{nnnnnnn.3} is satisfied. However, a positive answer to the
  conjecture \cite[Remark 5.20]{GGOR} would imply, via the \KZ-functor, that $F$ is isomorphic to
  the identity functor on the additive subcategory given by
  projective-injective objects. Since the corresponding Hecke algebra is symmetric
  (see \eg\cite[Lemma 5.10]{CIK}), the conjecture would follow from Proposition~\ref{cperfect}.\\

Independently of the Conjecture~\ref{conjCheredsym}, we can at least give a 
description of the corresponding Serre functor in
terms of partial coapproximation:

\begin{prop}
     Let $\op{H}(V,\gamma)$ be a rational Cherednik algebra with the corresponding
  category $\cO(V,\gamma)$. Let $Q$ be a basic projective-injective module in
  $\cO(V,\gamma)$. Then the Serre functor of $\cD^b(\cO(V,\gamma))$ is isomorphic
  to $\cL((\CA_Q)^2)$.
\end{prop}

\begin{proof}
  We only have to verify that we are in the situation of
  Theorem~\ref{Serrecoapprox}. We can find a projective-injective module
  $P_{KZ}$ representing the \KZ-functor (see \cite[Proposition
  5.21]{GGOR}). On the other hand the endomorphism ring of $P_{KZ}$ is
  isomorphic to the Hecke algebra (\cite[Theorem 5.15]{GGOR}), hence
  symmetric (\cite{CIK}). It is known that  $\cO(V,\gamma)$ has the double
  centraliser property with respect to $P_{KZ}$ (\cite[Theorem
  5.16]{GGOR}). Since the naive duality maps a basic projective module to a
  basic injective module, and $P_{KZ}$ to the corresponding
  $P_{KZ}$ (\cite[Proposition 5.21]{GGOR}), and $\cO(V,\gamma)^{opp}$ has again a double
  centraliser property with respect to $P_{KZ}$, we get that $A$ and $A^{opp}$
  have the double centraliser with respect to a basic projective-injective
  module. From \cite[Proposition 5.21]{GGOR} it follows that such a
  basic projective-injective module for $A$ is good. Hence, the assumptions 
  of Theorem~\ref{Serrecoapprox} are satisfied. The statement follows.    
\end{proof}

\section{Projective-injectives in the ca\-te\-go\-ry
  $\cO^{\p}(\mathfrak{sl}_n(\C))$} \label{s101}

In the following section we study more carefully projective-injective modules
in the parabolic category $\cO^{\p}$, especially, for the Lie algebra
$\mathfrak{sl}_n=\mathfrak{sl}_n(\C)$.  As already mentioned in the
introduction, one of the motivations to consider the category of
projective-injective modules in $\cO^{\p}$ is to find a precise
connection between Khovanov's categorification of the Jones polynomial 
(\cite{Khovanov}) and the categorification of the Jones polynomial via representation 
theory of the Lie algebra $\mathfrak{sl}_n$ (as proposed in \cite{BFK} and
proved in \cite{Sduke}). It might be possible to  pass directly 
from one model to the other by connecting the involved algebras directly, 
because the algebra, used by Khovanov in his categorification, is a quotient
of an algebra $A$ such that $A\MOD$ is equivalent to a certain block of
$\cO^{\p}$ for some $\mathfrak{sl}_n$ (\cite[page~494]{Braden}). Although, we have a very nice, 
more or less explicit, description of the algebra $A$ in question 
(\cite[Theorem~1.4.1]{Braden}), we are interested in more conceptual properties of 
the algebra. Several conjectures in this direction were formulated by Khovanov in
\cite{Kh1}. We want to simplify the problem by using the double centraliser
property. In this way, by using the Serre functor, we confirm
three conjectures of Khovanov: in Theorem~\ref{t5} we confirm \cite[Conjecture~4]{Kh1} 
concerning the centre of $A$, and in Theorem~\ref{tparab} we confirm that the
endomorphism algebra of a basic projective-injective module is symmetric, and
depends only on the chosen partition of $n$, not on the actually chosen
composition of $n$. (The last two conjectures were formulated in a private
communication). Furthermore, Theorem~\ref{tparab} supports \cite[Conjecture~3]{Kh1}.
\subsection{On a result of Irving}
Consider the classical triangular decomposition 
$\mathfrak{sl}_n=\mathfrak{n}_-\oplus\mathfrak{h}\oplus \mathfrak{n}_+$, 
where $\mathfrak{h}$ is the Cartan subalgebra of all diagonal
matrices (with zero trace) and $\mathfrak{n}_{\pm}$ denotes the subalgebra of 
all upper- and lower-triangular matrices respectively. Given a composition, $\mu=(\mu_1,\dots,\mu_k)$, of $n$ (\ie $\mu_1+\dots+\mu_k=n$), 
we have the corresponding Young subgroup $S_\mu=S_{\mu_1}\times S_{\mu_2}\times\dots
\times S_{\mu_k}$ of $S_n$, the latter being the Weyl group of $\mathfrak{sl}_n$.
Let $\mathcal{O}_0^{\mu}$ be the corresponding parabolic subcategory of 
$\cO_0(\mathfrak{sl}_n)$ (see Section~\ref{parabolic}). Recall that the simple objects in $\mathcal{O}_0$ are in bijection with the elements of $S_n$. 
We denote by $L(w)$ the 
simple module of highest weight $w(\rho)-\rho$, where $\rho$ denotes the half-sum of
the positive roots. 
For $w\in S_n$ we denote by $\Delta(w)$ and  $P(w)$ the corresponding
Verma and indecomposable projective module in $\mathcal{O}_0$.
Let $S^\mu$ be the set of shortest coset representatives of
$S_{\mu}\backslash S_n$ and let $w_{\mu}$ be the longest element in $S^\mu$.  
The simple objects in $\mathcal{O}_0^{\mu}$ are then the $L(w)$, where $w\in S^\mu$.
For $w\in S^\mu$ we denote by $L^{\mu}(w)$, $\Delta^{\mu}(w)$ and  $P^{\mu}(w)$ the 
corresponding simple, parabolic Verma and indecomposable projective module in 
$\mathcal{O}_0^{\mu}$  respectively. Note that $L^{\mu}(w)=L(w)$ for $w\in S^\mu$. 

For any $i\in\{1,\dots,n-1\}$ we denote by $\theta_i:\mathcal{O}_0\to
\mathcal{O}_0$ the translation functor through the $s_i$-wall (see \eg \cite[Section 3]{GJ}). 
This functor is exact, self-adjoint, and preserves $\mathcal{O}_0^{\mu}$. 
For $w\in S_n$ we denote by $\mathcal{L}(w)$ the left cell of the element $w$
(for a definition we refer to \cite{KLHecke}). Now we can give (for the
$\mathfrak{sl}_n$ case) an easier proof for the
following main result of \cite{Irself}:

\begin{theorem}\label{t7}
For any composition,  $\mu$, of $n$ the following conditions are equivalent:
\begin{enumerate}[(i)]
\item\label{co1} $P^{\mu}(w)$ is injective.
\item\label{co2} $w\in \mathcal{L}(w_{\mu})$
\item\label{co3} $L^{\mu}(w)$ occurs in the socle of some parabolic
Verma module $\Delta^{\mu}(w')$.
\end{enumerate}
\end{theorem}

\begin{proof}
Assume that \eqref{co1} is satisfied. 
Since $\mathcal{O}_0^{\mu}$ has a simple preserving duality,
if $P^{\mu}(w)$ is injective, it is a tilting module in the
highest weight category $\mathcal{O}_{0}^{\mu}$ and hence is self-dual. This means that its socle is $L^{\mu}(w)$, which must
coincide with the socle of some parabolic Verma module because  
$P^{\mu}(w)$ has a standard flag. This implies \eqref{co3}.

Assume that \eqref{co3} is satisfied. Since any parabolic Verma
module is a submodule of some tilting module, we get that
$L^{\mu}(w)$ occurs in the socle of some tilting module. 
By \cite{CI}, the tilting modules in $\mathcal{O}_{0}^{\mu}$ are exactly
direct summands of translations of $L^{\mu}(w_{\mu})$. From 
\cite[Proposition~4.3]{Irself} (this is an easy 
preparatory result) it follows that $w\in \mathcal{L}(w_{\mu})$, 
that is \eqref{co2} is satisfied.

Assume that \eqref{co2} is satisfied.
Since all $P^{\mu}(w')$, $w'\in \mathcal{L}(w_{\mu})$, can be
obtained from each other via translations through walls (this follows again 
from \cite[Proposition~4.3]{Irself}),
it is now left to show that there exists $w\in \mathcal{L}(w_{\mu})$
such that $P^{\mu}(w)$ is injective. Actually, since we already
know that \eqref{co1} implies \eqref{co2}, it is enough to show that
there exists some projective-injective module in $\mathcal{O}_{0}^{\mu}$.
But this one is obtained by translating any simple projective module
from the same weight lattice, which exists by \cite[3.1]{IS}
(this is again an easy result).
\end{proof}

\subsection{On Khovanov's conjectures}
According to Theorem~\ref{t7}, the modules $P^{\mu}(w)$,
$w\in \mathcal{L}(w_{\mu})$, constitute an exhaustive list of indecomposable
projective-injective modules in $\mathcal{O}_0^{\mu}$. Let
$P_{\mu}=\bigoplus_{w\in \mathcal{L}(w_{\mu})}P^{\mu}(w)$ be the basic
projective-injective module and set
$B_{\mu}=\mathrm{End}_{\mathcal{O}_0^{\mu}}\left(P_{\mu}\right)$.
As a consequence, we have the following result which confirms \cite[Conjecture~4]{Kh1}:

\begin{theorem}\label{t5}
Let $\mu$ be a composition of $n$.
\begin{enumerate}[(1)]
\item\label{t5.1} The Loewy lengths of all projective-injective modules in $\cO_0^\mu$
coincide. 
\item\label{t5.2}  $\cO_0^\mu$ satisfies the double centraliser property with respect to a
basic pro\-jec\-tive-\-in\-jec\-tive module, in particular, the restriction 
induces an isomorphism between the centres of $A^{\mu}_0$ and $B_{\mu}$.
\end{enumerate}
\end{theorem}

\begin{remark}
{\rm
Theorem~\ref{t5} is true for any semisimple complex Lie algebra $\mg$
(the proof is exactly the same if one replaces Theorem~\ref{t7} by the main result of 
\cite{Irself}).
}
\end{remark}

\begin{proof}
Let $P^{\mu}(w)$ be a projective-injective module.
From \cite[Proposition~4.3(ii)]{Irself} it follows that 
$P^{\mu}(w)$ is a direct summand of some translation of $P^{\mu}(w_{\mu})$.
Hence to prove the first statement it is enough to show that translations
through walls do not increase the Loewy length of projective-injective
modules. If $\theta_i(L(w))\neq 0$, then
\begin{displaymath}
\dim\HOM_{\cO}(\theta_iP^{\mu}(w),L(w))=
\dim\HOM_{\cO}(P^{\mu}(w),\theta_iL(w))=2
\end{displaymath}
by \cite[4.12(3),~4.13(3')]{Ja2}, which implies that $P^{\mu}(w)\oplus P^{\mu}(w)$ 
is a direct summand of $\theta_i(P^{\mu}(w))$. Comparing the lengths of the standard 
filtrations we even get $\theta_i(P^{\mu}(w))\cong P^{\mu}(w)\oplus P^{\mu}(w)$, 
in particular, such translations do not produce new projective-injective modules.

Now assume that $\theta_i(L(w))=0$.
The algebras $A_0$ and $A_0^{\mu}$, which correspond to $\mathcal{O}_0$ and
$\mathcal{O}_0^{\mu}$ are Koszul (\cite{Sperv}, \cite{BGS}, \cite{ErikKoszul}), 
in particular, they admit a canonical positive grading (the Koszul grading), 
which we fix. This allows us to consider graded versions of both $\mathcal{O}_0$ 
and $\mathcal{O}_0^{\mu}$ (see \cite{BGS}, \cite{Stgrad}).
In \cite{Stgrad} and \cite{BGS} it was shown that simple modules, Verma
modules, parabolic Verma and projective modules in $\mathcal{O}_0$ and 
$\mathcal{O}_0^{\mu}$ are gradable.  Their graded lifts are unique up to isomorphism 
and grading shift, therefore we call a lift {\it standard} if the head is concentrated in
degree zero. In \cite{Stgrad} it was shown, that the functors $\theta_i$ 
(as endofunctors of $\cO_0$) are gradable as well. 
We denote by $\tilde{\theta}_i$ the standard graded lift
of $\theta_i$ (\ie $\tilde{\theta}_i$, applied to a simple module concentrated
in degree $0$ has socle concentrated in degree $1$). Since
$\theta_i$ preserves $\mathcal{O}_0^{\mu}$, the functor
$\tilde{\theta}_i$ restricts to a graded lift of $\theta_i$ on
$\mathcal{O}_0^{\mu}$.

Let $P^{\op{gr}}$  be the standard graded lift 
of $P^{\mu}(w)$. Since it has both simple top and simple socle, the radical-, socle- 
and graded filtrations of $P^{\op{gr}}$ coincide by \cite[Proposition~2.4.1]{BGS}. 
In particular, $P^{\op{gr}}$ has a unique component of maximal and a unique 
component of minimal degree. On the other hand, $\tilde{\theta}_i(L)$ is concentrated
in the degrees $-1,0,1$ for any simple module $L$, concentrated in degree $0$
(\cite[Theorem~5.1]{Stgrad}). This implies that the
length of the graded filtration of $\tilde{\theta}_i(P^{\op{gr}})$ can not 
exceed the length of the graded filtration of $P^{\op{gr}}$. Hence,
the Loewy length of $\theta_i(P^{\mu}(w))$ does not exceed that of $P^{\mu}(w)$ and
the statement \eqref{t5.1} follows.

The double centraliser property follows from Corollary~\ref{nice} and the main 
result of \cite{Irself} (as formulated in Theorem~\ref{t7}). 
For an algebra, $A$, we denote its centre by $Z(A)$. From the 
double centraliser property we have 
\begin{displaymath}
A_0^{\mu}=\END_{B_{\mu}}(P_{\mu}).
\end{displaymath}
If $x\in Z(A_0^{\mu})$, then $xa=ax$ for all $a\in A_0^{\mu}$ and hence
$x\in B_{\mu}$. On the other hand, every element of $A_0^{\mu}$ commutes with
each element in $B_{\mu}$ by definition. Hence $x\in Z(B_{\mu})$. This
implies $Z(A_0^{\mu})\hookrightarrow Z(B_{\mu})$. Because of the left-right
symmetry of the double centraliser
we finally get $Z(A_0^{\mu})= Z(B_{\mu})$.
This completes the proof.
\end{proof} 

We formulate now the main result of this section: 
\begin{theorem}\label{tparab}
\begin{enumerate}[(1)]
\item\label{tparab.1} Let $\mu$ be a composition of $n$. Then the algebra $B_{\mu}$
is symmetric.
\item\label{tparab.2} Let $\mu$ and $\nu$ be two compositions of $n$, which give rise 
to the same partition of $n$. Then $B_{\mu}\cong B_{\nu}$.
\end{enumerate}
\end{theorem}

\begin{proof}
The statement~\eqref{tparab.1} is just a special case of
Theorem~\ref{symmetry}. 

Let us now prove \eqref{tparab.2}.
Without loss of generality we may assume $\mu=(\mu_1,\dots,\mu_k)$ and 
\begin{displaymath}
\nu=(\mu_1,\dots,\mu_{l-1},\mu_{l+1},
\mu_{l},\mu_{l+2},\dots,\mu_k)
\end{displaymath}
for some $l\in\{1,\dots,k-1\}$. Moreover, we assume $\mu_l>\mu_{l+1}$.
 
For any composition, $\tau$, of $n$ Irving and Shelton constructed in
\cite[3.1]{IS} a special weight, $\lambda(\tau)$, with the following
property: the simple highest weight module $L(\lambda(\tau)-\rho)$ 
with the highest weight $\lambda(\tau)-\rho$ is the only simple module in its
block of $\cO^{\tau}$. From the definition in \cite[3.1]{IS} it follows 
immediately that, if $\tau$ and  $\tau'$ are two compositions of $n$ 
which give rise to the same partition, then
$\lambda(\tau)$ and $\lambda(\tau')$ are in the same $S_n$-orbit,
in particular, $L(\lambda(\tau)-\rho)$ and 
$L(\lambda(\tau')-\rho)$ belong to the same block of $\cO$. 

Now we apply this to the case $\tau=\nu$, $\tau'=\mu$.
Let $\cO_{\xi}$ be the common block (of $\cO$) for $L(\lambda(\mu)-\rho)$
and $L(\lambda(\nu)-\rho)$. Although $L(\lambda(\mu)-\rho)$
and $L(\lambda(\nu)-\rho)$ are in the same block of $\cO$, we have that the 
parabolic categories $\cO_{\xi}^{\mu}$
and $\cO_{\xi}^{\nu}$ are semi-simple containing only one simple object each.
Obviously, they are equivalent. However, we would like to construct a
functor on $\cO$, which gives rise to an equivalence 
between these categories, and, additionally, commutes with tensoring
with finite-dimensional $\mathfrak{sl}_n$-modules.

To proceed we will need some general notation. For any transposition
$s=s_i=(i,i+1)$ in $S_n$ we set $\cO^s=\cO^{\beta}$, where 
$\beta=(\beta_1,\dots,\beta_{n-1})$ is the composition of $n$ such that 
$\beta_i=2$ and $\beta_j=1$ for all $j\neq i$. Denote by 
$\mathfrak{i}_s:\cO^s\hookrightarrow\cO$ the inclusion functor, and
by $\mathrm{Z}_s:\cO\to \cO^s$ the left adjoint to $\mathfrak{i}_s$, which is
the {\em Zuckerman functor}, associated to $s$. Then 
$\op{d}\mathrm{Z}_{s}\op{d}:\cO^s\hookrightarrow\cO$ is the right adjoint to 
$\mathfrak{i}_s$. It is known that 
\begin{equation}\label{eqew}
\mathcal{L}_i\mathrm{Z}_{s}=0\text{  for all } 
i\geq 3,\quad\text{ and }\quad
\mathcal{L}\mathrm{Z}_{s}\cong \op{d}\mathcal{L}\mathrm{Z}_{s}\op{d}[2],
\end{equation}
see \cite{EW} or \cite[Proposition~3]{BFK} and also Lemma~\ref{Zucker}. Finally, we denote by
$\mathfrak{i}_{\mu}:\cO^{\mu}\hookrightarrow\cO$ and
$\mathfrak{i}_{\nu}:\cO^{\nu}\hookrightarrow\cO$ the inclusion functors and
by $\mathrm{Z}_{\mu}:\cO\to \cO^{\mu}$ and $\mathrm{Z}_{\nu}:\cO\to \cO^{\nu}$
the corresponding left adjoint Zuckerman functors.

Let $m'=m'_l=\mu_1+\dots+\mu_{l-1}$ and $m=m_l=m'+\mu_{l+1}$, and set 
$r=\mu_l-\mu_{l+1}$. Consider the following element in $S_n$:
\begin{displaymath}
w=(s_{m'+r+1}\dots s_{m+r})\dots(s_{m'+2}\dots s_{m+1}s_{m'+1})
(s_{m'+1}s_{m'+2}\dots s_{m-1}s_m).
\end{displaymath}
(for example if $\mu=(2,4)$ then $\nu=(4,2)$ and
we have $w=(s_2s_3)(s_1s_2)$).
For simplicity we write the above product in the form $w=t_1\dots t_{p}$
(thus $p=r\mu_{l+1}$, $t_1=s_{m'+r+1}$ and so on). For $i=1,\dots,p$ we set
$w_i=t_i\dots t_{p}$, and $w_{p+1}=e$. The element $w$ is constructed
such that $w(\lambda(\nu))=\lambda(\mu)$ and $w_i(\lambda(\nu))>w_{i+1}(\lambda(\nu))$
for all $i=1,\dots,p$.  Define
\begin{displaymath}
\begin{array}{rclcl}
F&=&\mathcal{L}\mathrm{Z}_{\mu}\mathfrak{i}_{t_1}\mathcal{L}\mathrm{Z}_{t_1}
\mathfrak{i}_{t_2}\mathcal{L}\mathrm{Z}_{t_2}\dots
\mathfrak{i}_{t_p}\mathcal{L}\mathrm{Z}_{t_p}\mathfrak{i}_{\nu}[-l(w)]&:&
\mathcal{D}^b(\cO^{\nu})\to \mathcal{D}^b(\cO^{\mu}), \\
G&=&\mathcal{L}\mathrm{Z}_{\nu}\mathfrak{i}_{t_p}\mathcal{L}\mathrm{Z}_{t_p}
\mathfrak{i}_{t_{p-1}}\mathcal{L}\mathrm{Z}_{t_{p-1}}\dots
\mathfrak{i}_{t_1}\mathcal{L}\mathrm{Z}_{t_1}\mathfrak{i}_{\mu}[-l(w)]&:&
\mathcal{D}^b(\cO^{\mu})\to \mathcal{D}^b(\cO^{\nu}).
\end{array}
\end{displaymath}
Both, $F$ and $G$, commute with tensoring with finite-dimensional 
$\mathfrak{sl}_n$-modules (\cite[Proposition~2.2 and Proposition~3.7]{EW}, see also 
\cite[Proposition~3]{BFK}). Further, $F$ is both left and right adjoint to
$G$ by \eqref{eqew}.\\

{\it Claim: The functors $F$ and $G$ define, via restriction, mutually inverse 
equivalences $F:\cO_{\xi}^{\nu}\to \cO_{\xi}^{\mu}$ and 
$G:\cO_{\xi}^{\mu}\to \cO_{\xi}^{\nu}$.}

\begin{proof}[Proof of the Claim.]
Since $F$ and $G$ are adjoint to each other and both $\cO_{\xi}^{\mu}$ and
$\cO_{\xi}^{\nu}$ are semi-simple, it is enough to show that
\begin{displaymath}
F(L(\lambda(\nu)-\rho))=L(\lambda(\mu)-\rho).
\end{displaymath}

To prove this we first note that for a simple reflection, $s$, and for
a dominant integral weight, $\lambda$,  we have
\begin{equation}\label{zucker}
\mathcal{L}\mathrm{Z}_{s}(L(x\cdot\lambda))=
\begin{cases}
L(x\cdot\lambda)\oplus L(x\cdot\lambda)[2], & sx\cdot\lambda< x\cdot\lambda, \\
\bigoplus_{y:y>sx} L(y\cdot\lambda)^{a_y}[1], & sx\cdot\lambda= x\cdot\lambda,\\
L(sx\cdot\lambda)[1]\oplus \bigoplus_{y:y>sx} L(y\cdot\lambda)^{a_y}[1], & \mathrm{otherwise},
\end{cases}
\end{equation}
where $a_y\in\{0,1,\dots\}$. To see this, let $\mathrm{T}_s:\cO\to\cO$ be the twisting
functor, associated with $s$ (as in Subsection~\ref{secappl.1}).  
In \cite[Theorem~4]{KM} (see also \cite[Proposition~2.3]{MazStr}) it is shown that
there exists a natural transformation, $\mathrm{can}_s:\mathrm{T}_s\to\mathrm{ID}$, 
non-vanishing on Verma modules. In \cite[Theorem~1(3)]{MazStr} it is proved that 
the kernel of $\mathrm{can}_s$ is isomorphic to $\mathcal{L}_1\mathrm{Z}_s$.
Now \eqref{zucker} follows from the Kazhdan-Lusztig conjectures, see
\cite[Theorem~6.3 and Theorem~7.8]{AS} for details. (Note that the assumption for
$L'$ to be $s$-finite is missing in the formulation of \cite[Theorem~6.3(3)]{AS}.)

For $i=1,\dots,p$ set
\begin{displaymath}
F_i=\mathcal{L}\mathrm{Z}_{t_i}
\mathfrak{i}_{t_2}\mathcal{L}\mathrm{Z}_{t_2}\dots
\mathfrak{i}_{t_p}\mathcal{L}\mathrm{Z}_{t_p}\mathfrak{i}_{\nu}[-(p-i+1)]:
\mathcal{D}^b(\cO^{\nu})\to \mathcal{D}^b(\cO^{t_i}).
\end{displaymath}
From \eqref{zucker} it follows by induction that
\begin{displaymath}
F_iL(\lambda(\tau)-\rho)=L(w_{i}(\lambda(\tau))-\rho)
\oplus \bigoplus_{y:y>w_{i}} L(y\cdot\lambda(\tau))^{a_y^{(i)}}[b_y^{(i)}],
\end{displaymath}
where $a_y^{(i)}\in\{0,1,\dots\}$, $b_y^{(i)}\in\{1,\dots\}$ for all $i$. 
In particular, we have
\begin{displaymath}
\mathcal{L}\mathrm{Z}_{\mu}\mathfrak{i}_{t_1}F_1L(\lambda(\tau)-\rho)=
L(\lambda(\mu)-\rho)
\end{displaymath}
since $\cO_{\xi}^{\mu}$ is semi-simple 
and the claim follows.
\end{proof}

Let us now sum up what we know. We have an adjoint  pair, $(F,G)$, of
functors between $\mathcal{D}^b(\cO^{\nu})$ and $\mathcal{D}^b(\cO^{\mu})$, 
which commute with tensoring with finite-dimensional $\mathfrak{sl}_n$-modules 
and induce mutually inverse equivalences, when restricted to
$\cO_{\xi}^{\nu}$ and $\cO_{\xi}^{\mu}$. On the other hand, there is a
finite dimensional $\mathfrak{sl}_n$-module, $E$, such that 
$E\otimes L(\lambda(\nu)-\rho)$ contains $P_{\nu}$ as a direct summand,
and $E\otimes L(\lambda(\mu)-\rho)$ contains $P_{\mu}$ as a direct summand
(this follows from Theorem~\ref{t5}, for the explicit statement see
\cite[Proposition~4.3(ii)]{Irself}). Therefore the adjunction morphisms
$FG\to\ID$ and $\ID\to GF$ are isomorphisms when evaluated at
$P_{\mu}$ and $P_{\nu}$ respectively. Hence $F$ and $G$ define mutually
inverse equivalences between the corresponding additive categories of
projective-injective modules.
This completes the proof of Theorem~\ref{tparab}.
\end{proof}

We have the following direct consequence, a part of which was also obtained 
in \cite[Section~6]{Kh1} by establishing a derived equivalence between
$\cO_{\xi}^{\mu}$ and $\cO_{\xi}^{\nu}$ using a geometric argument:

\begin{cor}\label{khovcenter}
The centres of $B_{\mu}$, $B_{\nu}$, $\cO_{\xi}^{\mu}$ and 
$\cO_{\xi}^{\nu}$ are all isomorphic.
\end{cor}

\begin{remark}\label{bred}
{\rm 
Since the Kazhdan-Lusztig left cell modules for the Iwahori-Hecke algebra
$\mathcal{H}_n$ of the symmetric group are exactly the irreducible modules, 
Theorem~\ref{t7} can be used to ``categorify'' these irreducible modules: 
Let $\la$ be a partition of $n$. Consider the abelian category of
modules, admitting a $2$-step presentation by projective-injective 
modules in the parabolic category $\cO$ for $\mathfrak{sl}_n$, associated
with $\lambda$. This abelian category is invariant under the action of
translations through walls. The action of these translation functors 
gives rise to a categorification of the Specht module $S_{\lambda}$ for the symmetric group $S_n$.
The graded version of the above result (in the sense of \cite{Stgrad})
gives rise to a categorification of the Specht module $S_{\lambda}$ for 
$\mathcal{H}_n$. The details will appear in \cite{KMS}.
}
\end{remark}

\bibliography{ref}

\noindent
Volodymyr Mazorchuk, Department of Mathematics, Uppsala University,
Box 480, 751 06, Uppsala, SWEDEN,\\ 
e-mail: {\tt mazor\symbol{64}math.uu.se},
web: {``http://www.math.uu.se/$\tilde{\hspace{1mm}}$mazor/''}.
\vspace{0.2cm}

\noindent
Catharina Stroppel, Department of Mathematics,
University of Glasgow, University Gardens,
Glasgow G12 8QW, UK,\\
e-mail: {\tt cs\symbol{64}maths.gla.ac.uk}
\vspace{0.2cm}

\end{document}